\documentclass[12 pt]{article}
\usepackage{amssymb,amsmath,amsfonts,amsthm}
\newcommand\Cov{{\mathrm {Cov}}}

\newcommand\bkR{{\mathbb {R}}}

\newcommand\bkE{{\mathbb {E}}}

\newtheorem {Lemma}{Lemma}[section]
\newtheorem {Theorem}{Theorem}[section]
\newtheorem {Proposition}{Proposition}[section]
\newtheorem {Corollary}{Corollary}[section]

\theoremstyle{definition}
\newtheorem{Definition}{Definition}[section]
\newtheorem{Notation}{Notation}[section]
\newtheorem{Remark}{Remark}[section]

\newcommand\I{{ 1\hspace{-1,2mm}{\mathrm I}}}

\newcommand\no{\noindent}
\newcommand\ssk{\smallskip}
\newcommand\msk{\medskip}

\newcommand\beq{\begin{equation}}
\newcommand\eeq{\end{equation}}

\def\ti{\tilde}\def\no{\noindent}
\def\msk{\medskip}\def\ssk{\smallskip}

\def\X{(X_i)_{i\in{\mathbb N}}}
\def\XZ{(X_i)_{i\in{\mathbb Z}}}
\def\Y{(Y_i)_{i\in{\mathbb N}}}

\begin{document}
\begin{center} {\bf \Large Rates of convergence for minimal distances in
the central limit theorem under projective criteria}\vskip15pt

J\'er\^ome Dedecker $^{a}$, Florence Merlev\`{e}de $^{b}$ {\it
and\/} Emmanuel Rio   $^{c}$
\end{center}
\vskip10pt $^a$ Universit\'e Paris 6, LSTA, 175 rue
du Chevaleret, 75013 Paris, FRANCE.\\ E-mail: dedecker@ccr.jussieu.fr\\ \\
$^b$ Universit\'e Paris 6, LPMA and C.N.R.S UMR 7599, 175 rue
du Chevaleret, 75013 Paris, FRANCE. E-mail: merleve@ccr.jussieu.fr\\ \\
$^c$ Universit\'e de Versailles, Laboratoire de math\'ematiques, UMR
8100 CNRS, B\^atiment Fermat, 45 Avenue des Etats-Unis, 78035
Versailles, FRANCE. E-mail: rio@math.uvsq.fr \vskip10pt


\begin{center}
{\bf Abstract}\vskip10pt
\end{center}

In this paper, we give estimates of ideal or minimal distances
between the distribution of the normalized partial sum and the
limiting Gaussian distribution for stationary martingale difference sequences
or stationary sequences satisfying projective criteria. Applications
to functions of linear processes and  to functions of expanding maps
of the interval are given.
\section{Introduction and Notations}
\setcounter{equation}{0} Let $X_1, X_2, \ldots $ be a strictly
stationary sequence of real-valued random variables (r.v.) with mean
zero and finite variance. Set $S_n = X_1 + X_2 + \cdots + X_n$. By
$P_{n^{-1/2}S_n}$ we denote the law of $n^{-1/2}S_{n}$ and by
$G_{\sigma^2}$  the normal distribution $ N(0, \sigma^2)$. In this
paper, we shall give  quantitative estimates of the approximation of
$P_{n^{-1/2}S_n}$ by $G_{\sigma^2}$ in terms of minimal or ideal
metrics.

Let ${\mathcal L}(\mu, \nu)$ be the set of the probability laws on
$\mathbb R^2$ with  marginals $\mu$ and $\nu$. Let us consider the
following minimal distances (sometimes called Wasserstein distances
of order $r$)
$$
W_r(\mu, \nu) =  \left\{
\begin{array}{ll} \displaystyle \inf \Big \{ \int |x-y|^r P(dx, dy)  :  P
\in {\mathcal L}(\mu, \nu) \Big \} & \text{if $0< r < 1$} \\
\displaystyle \inf \Big \{ \Big (\int |x-y|^r P(dx, dy) \Big )^{1/r}
: P \in {\mathcal L}(\mu, \nu) \Big \}  & \text{if $r \geq 1$} \, .
\end{array}
\right.
$$
It is well known that for two probability measures $\mu$ and $\nu$
on $\mathbb R$ with respective distributions functions (d.f.) $F$
and $G$,
\begin{equation} \label{def2wasser}
W_r(\mu, \nu) =  \Big ( \int_0^1 |F^{-1}(u) - G^{-1}(u) |^r du\Big
)^{1/r} \, \text{ for any $r \geq 1$.}
\end{equation}

We consider also the following ideal distances of order $r$
(Zolotarev distances of order $r$). For two probability measures
$\mu$ and $\nu$, and $r$ a positive real, let
$$
\zeta_{r} ( \mu , \nu )   =  \sup \Big \{ \int fd\mu - \int f d \nu
: f \in \Lambda_r \Big \} \, ,
$$
where $\Lambda_r$ is defined as follows: denoting by $l$ the natural
integer such that $ l < r \leq l +1$, $\Lambda_r$ is the class of
real functions $f$ which are $l$-times continuously differentiable
and such that \beq \label{defgamrrio} |f^{(l)}(x) - f^{(l)}(y) |
\leq | x - y |^{r-l} \ \hbox{ for any } (x,y) \in \bkR \times \bkR
\, . \eeq It follows from the  Kantorovich-Rubinstein theorem (1958)
that for any $0 < r \leq 1$, \beq \label{KRdual} W_r(\mu, \nu) =
\zeta_{r} ( \mu , \nu ) \, . \eeq For probability laws on the real
line, Rio (1998) proved that  for any $r
> 1$,
\beq \label{riolink} W_r(\mu, \nu) \leq c_r \big ( \zeta_{r} ( \mu ,
\nu ) \big )^{1/r} \, , \eeq where $c_r$ is a constant depending
only on $r$.

For independent random variables, Ibragimov (1966) established that
if $X_1 \in {\mathbb L^p}$ for $p \in ]2,3]$, then $ W_1 (
P_{n^{-1/2}S_n},G_{\sigma^2}) = O(n^{1 - p/2})$ (see his Theorem
4.3). Still in the case of independent r.v.'s, Zolotarev (1976)
obtained the following upper bound for the ideal distance: if $X_1
\in {\mathbb L^p}$ for $p \in ]2, 3]$, then $\zeta_p
(P_{n^{-1/2}S_n},G_{\sigma^2}) = O(n^{1 - p/2})$. From
(\ref{riolink}), the result of Zolotarev entails that, for $p \in
]2, 3]$, $W_p ( P_{n^{-1/2}S_n},G_{\sigma^2}) = O(n^{1/p - 1/2})$
(which was obtained by Sakhanenko (1985) for any $p > 2$). From
(\ref{def2wasser}) and H\"older's inequality, we easily get that for
independent random variables in ${\mathbb L}^p$ with $p \in ]2, 3]$,
\beq \label{boundWrind} W_r ( P_{n^{-1/2}S_n},G_{\sigma^2}) =
 O(n^{-(p-2)/2 r}) \quad \text{for any $1 \leq r \leq p$}.
\eeq

In this paper, we are interested in extensions of (\ref{boundWrind})
to sequences of dependent random variables.  More precisely, for
$X_1 \in {\mathbb L}^p$ and $p$ in $]2, 3]$ we shall give ${\mathbb
L}^p$-projective criteria under which:
 for  $r \in [p-2,p]$ and $(r,p)\neq (1,3)$, \beq
\label{bound2Wrind} W_r ( P_{n^{-1/2}S_n},G_{\sigma^2}) =
 O(n^{-(p-2)/2\max(1,r)}) \, .
\eeq
As we shall see in Remark \ref{rmkcompprok},
(\ref{bound2Wrind}) applied to $r=p-2$ provides the rate of
convergence $O(n^{-\frac{p-2}{2(p-1)}})$ in the Berry-Esseen
theorem.

When $(r,p) = (1,3)$, Dedecker and Rio (2007)  obtained that
$W_1 ( P_{n^{-1/2}S_n},G_{\sigma^2}) = O(n^{-1/2})$
for stationary sequences of random variables in ${\mathbb L}^3$ satisfying
${\mathbb L}^1$ projective criteria or weak dependence
 assumptions (a similar result was obtained by P\`ene (2005) in the case where the variables are bounded). In
 this particular case our approach provides a new
 criterion under which $W_1 ( P_{n^{-1/2}S_n},G_{\sigma^2}) =
 O(n^{-1/2}\log n )$.

Our paper is organized as follows. In Section \ref{SectionMart}, we
give projective conditions for stationary martingales differences
sequences to satisfy (\ref{bound2Wrind}) in the  case $(r,p)\neq
(1,3)$. To be more precise, let $(X_i)_{i \in {\mathbb Z}}$ be a
stationary sequence of martingale differences with respect to some
$\sigma$-algebras $({\mathcal F}_i)_{i \in {\mathbb Z}}$ (see
Section \ref{nota} below for the definition of $({\mathcal F}_i)_{i
\in {\mathbb Z}}$). As a consequence of our Theorem \ref{Thmart}, we
obtain that if $(X_i)_{i \in {\mathbb Z}}$ is in ${\mathbb L}^p$
with $p \in ]2 , 3 ]$ and satisfies \beq \label{condmartintro}
\sum_{n=1}^{\infty}\frac{1}{n^{2-p/2}}\Big \Vert \bkE \Big
(\frac{S_n^2}{n} \Big | {\mathcal F}_0 \Big ) - \sigma^2 \Big
\Vert_{p/2} < \infty \, ,\eeq then the upper bound
(\ref{bound2Wrind}) holds provided that $(r,p)\neq (1,3)$. In the
case $r=1$ and $p=3$, we obtain the upper bound
$W_1 (
P_{n^{-1/2}S_n},G_{\sigma^2}) =
 O(n^{-1/2}\log n )$.

  In Section
 \ref{Sectionstatseq}, starting from the
 coboundary decomposition going back to Gordin (1969), and using the
 results of Section \ref{SectionMart}, we obtain ${\mathbb L}^p$-projective criteria ensuring
 (\ref{bound2Wrind}) (if $(r,p)\neq (1,3)$). For instance, if $(X_i)_{i \in {\mathbb
Z}}$ is a stationary sequence of ${\mathbb L}^p$ random variables
adapted to  $({\mathcal F}_i)_{i \in {\mathbb Z}}$, we obtain
(\ref{bound2Wrind}) for any $p\in ]2,3[$ and any $r \in [p-2,p]$
provided that (\ref{condmartintro}) holds and the series $\bkE (S_n
| {\mathcal F}_0) $ converge in ${\mathbb L}^p$. In the case where
$p=3$, this last condition has to be strengthened. Our approach
makes also possible to treat the case of non-adapted sequences.

 Section \ref{SectionAppli} is devoted to applications. In particular, we
give sufficient conditions for some functions of Harris recurrent
Markov chains and for functions of linear processes to satisfy the
bound (\ref{bound2Wrind}) in the case $(r,p)\neq (1,3)$ and the rate
$O(n^{-1/2}\log n )$ when $r=1$ and $p=3$. Since projective
criteria are verified under weak dependence assumptions, we give an
application to functions of $\phi$-dependent sequences in the sense
of Dedecker and Prieur (2007). These conditions apply to unbounded
functions of  uniformly expanding maps.

\medskip

\subsection{Preliminary notations}\label{nota} Throughout the paper, $Y$ is a
$N(0,1)$-distributed random variable.  We shall also use the
following notations. Let $(\Omega,{\cal A}, {\mathbb P} )$ be a
probability space, and $T:\Omega \mapsto \Omega$ be
 a bijective bimeasurable transformation preserving the probability ${\mathbb P} $.
For a $\sigma$-algebra ${\cal F}_0 $ satisfying ${\cal F}_0
\subseteq T^{-1 }({\cal F}_0)$, we define the nondecreasing
filtration $({\cal F}_i)_{i \in {\mathbb Z}}$ by ${\cal F}_i =T^{-i
}({\cal F}_0)$. Let ${\cal {F}}_{-\infty} = \bigcap_{k \in {\mathbb
Z}} {\cal {F}}_{k}$ and ${\cal {F}}_{\infty} = \bigvee_{k \in
{\mathbb Z}} {\cal {F}}_{k}$. We shall denote sometimes by ${\mathbb
E}_i$ the conditional expectation with respect to ${\mathcal F}_i$.
Let $X_0$ be a zero mean random variable with finite variance, and
define the stationary sequence  $(X_i)_{i \in \mathbb Z}$ by
 $X_i = X_0 \circ T^i$.\\

\section{Stationary sequences of martingale differences.}
\label{SectionMart}

\setcounter{equation}{0} In this section we give  bounds
for the ideal distance of order $r$ in the central limit theorem for
stationary martingale differences sequences $(X_i)_{i \in {\mathbb
Z}}$ under projective conditions.

\begin{Notation}
For any $p > 2$, define the envelope norm $\Vert \, . \, \Vert_{1,
\Phi, p}$ by
$$
\Vert X \Vert_{1, \Phi, p}  = \int_0^1 (1 \vee \Phi^{-1} (1-u/2)
)^{p-2} Q_X(u) du
$$
where $Q_X$ denotes the quantile function of $|X|$, and $\Phi$
denotes the d.f. of the $N(0, 1)$ law.
\end{Notation}
\begin{Theorem} \label{Thmart} Let $(X_i)_{i \in \mathbb Z}$ be a stationary
martingale differences sequence with respect to $({\cal F}_i)_{i
\in {\mathbb Z}}$. Let $\sigma$ denote the standard deviation of
$X_0$. Let $p \in ]2,3]$. Assume that $\bkE |X_0 |^p < \infty$ and
that
 \beq \label{C1}
\sum_{n=1}^{\infty}\frac{1}{n^{2-p/2}}\Big \Vert \bkE \Big
(\frac{S_n^2}{n} \Big | {\mathcal F}_0 \Big ) - \sigma^2 \Big
\Vert_{1, \Phi, p} < \infty \, ,\eeq and \beq \label{C2}
\sum_{n=1}^{\infty}\frac{1}{n^{ 2/p}}\Big \Vert \bkE \Big
(\frac{S_n^2}{n} \Big| {\mathcal F}_0 \Big ) - \sigma^2
\Big\Vert_{p/2} < \infty \, .\eeq Then,  for any $r \in [p-2, p]$
with $(r,p) \neq (1,3)$, $ \zeta_{r}( P_{n^{-1/2}S_n},G_{\sigma^2}
)=
 O(n^{1-p/2})$, and for $p=3$
 $\zeta_{1}( P_{n^{-1/2}S_n},G_{\sigma^2} )=
O( n^{-1/2}\log n )$.
\end{Theorem}
\begin{Remark} \label{rmkr<p-2}
Under the assumptions of Theorem \ref{Thmart}, $ \zeta_{r}(
P_{n^{-1/2}S_n},G_{\sigma^2} )= O (n^{-r/2})$ if $r<p-2$. Indeed,
let $p'=r+2$. Since $p' < p$, if the conditions (\ref{C1}) and
(\ref{C2}) are satisfied for $p$, they also hold for $p'$. Hence
Theorem \ref{Thmart} applies with $p'$.
\end{Remark}
From (\ref{KRdual}) and (\ref{riolink}), the following result
holds for the Wasserstein distances of order $r$.
\begin{Corollary} Under the conditions of Theorem \ref{Thmart},
$W_r ( P_{n^{-1/2}S_n},G_{\sigma^2}) =  O(n^{-(p-2)/ 2 \max (1,
r)})$ for any $r$ in $[p-2,p]$, provided that $(r,p) \not= (1,3)$.
\end{Corollary}
\begin{Remark}
For $p$ in $]2,3]$,  $W_p ( P_{n^{-1/2}S_n},G_{\sigma^2}) =
O(n^{-(2-p)/2p})$. This bound was obtained by Sakhanenko (1985) in
the independent case. For $p<3$, we have $W_1 (
P_{n^{-1/2}S_n},G_{\sigma^2}) = O(n^{1 - p/2})$. This bound was
obtained by Ibragimov (1966) in the independent case.
\end{Remark}
\begin{Remark} \label{rmkcompprok}
Let $\Pi_n$ be the Prokhorov distance between the law of
$n^{-1/2}S_n$ and the normal distribution $ N(0, \sigma^2)$.
From Markov's inequality,
$$
\Pi_n  \leq ( W_r(P_{n^{-1/2}S_n},G_{\sigma^2}) )^{1/(r+1)} \text{
for any $0 < r\leq 1$} \, .
$$
Taking $r=p-2$, it  follows  that  under the assumptions of Theorem
\ref{Thmart}, \beq \label{majprokh} \Pi_n =
O(n^{-\frac{p-2}{2(p-1)}}) \ \ \text{if $p < 3$ and} \ \Pi_n= O(
n^{-1/4}\sqrt{\log n}) \ \ \text{ if $p=3$.} \eeq For $p$ in
$]2,4]$, under (\ref{C2}), we have that $\Vert \sum_{i=1}^n \bkE
(X_i^2 - \sigma^2 |{\cal F}_{i-1}) \Vert_{p/2}=O(n^{2/p} ) $ (apply
Theorem 2 in Wu and Zhao (2006)). Applying then the result in Heyde
and Brown (1970), we get that if $\XZ$ is a stationary martingale
difference sequence in ${\mathbb L}^p$ such that  (\ref{C2}) is
satisfied then
$$
\Vert F_n -\Phi_\sigma \Vert_{\infty}= O \big (
n^{-\frac{p-2}{2(p+1)}} \big ) \, .
$$
where $F_n$ is the distribution function of $n^{-1/2} S_n$ and
$\Phi_\sigma$ is the d.f. of $G_{\sigma^2}$. Now
$$
 \Vert F_n -\Phi_\sigma \Vert_{\infty} \leq
\big(1 + \sigma^{-1}  (2 \pi)^{-1/2} \big)  \Pi_n \, .
$$
 Consequently the bounds obtained in (\ref{majprokh}) improve the one given in Heyde and
Brown (1970), provided that (\ref{C1}) holds.
\end{Remark}
\begin{Remark}
Notice that if $\XZ$ is a stationary martingale difference sequence
in ${\mathbb L}^3$ such that  ${\mathbb E} (X_0^2) = \sigma^2$ and
\begin{equation}\label{weakJan}
\sum_{k>0} k^{-1/2} \Vert \bkE (X_k^2|{\mathcal F}_0 )-
\sigma^2 \Vert_{3/2} < \infty ,
\end{equation}
then the conditions (\ref{C1}) and (\ref{C2}) hold for $p=3$.
Consequently, if  (\ref{weakJan}) holds, then Remark
\ref{rmkcompprok} gives $\Vert F_n -\Phi_\sigma \Vert_{\infty}= O
\big ( n^{-1/4} \sqrt{\log n} \big )$. This result has to be
compared with Theorem 6 in Jan (2001), which states that $\Vert F_n
-\Phi_\sigma \Vert_{\infty} = O (n^{-1/4}  )\,$  if   $\, \sum_{k>0}
\Vert \bkE (X_k^2|{\mathcal F}_0 )- \sigma^2 \Vert_{3/2} < \infty $.
\end{Remark}
\begin{Remark} \label{rmkequivalence}
Notice that if $(X_i)_{i \in \mathbb Z}$ is a stationary martingale
differences sequence, then the conditions (\ref{C1}) and (\ref{C2})
are  respectively equivalent to
$$
\sum_{j\geq 0}2^{j(p/2-1)} \Vert 2^{-j}\, \bkE  (S_{2^j}^2  |
{\mathcal F}_0) - \sigma^2  \Vert_{ 1 , \Phi , p} < \infty, \
\text{and} \  \sum_{j\geq 0}2^{ j(1-2/p)}\Vert 2^{-j}\, \bkE
(S_{2^j}^2  | {\mathcal F}_0) - \sigma^2  \Vert_{p/2} < \infty \, .
$$
To see this,  let $
A_n = \Vert \bkE (S_n^2 | {\mathcal F}_0 \big ) - \bkE
(S_n^2)\Vert_{1, \Phi, p}$ and   $B_n = \Vert \bkE (S_n^2 |
{\mathcal F}_0 \big ) - \bkE (S_n^2)\Vert_{p/2}$. We first show that
$A_n$ and $B_n$ are subadditive sequences. Indeed, by the martingale
property and the stationarity of the sequence,  for all positive $i$ and $j$
\begin{eqnarray*}
A_{i+j} &=& \Vert \bkE ( S_i^2 + (S_{i+j}-S_i)^2 | {\mathcal F}_0
\big ) - \bkE (S_i^2 + ( S_{i+j}-S_i)^2 ) \Vert_{1, \Phi, p}\\&\leq
& A_i + \Vert \bkE \big( ( S_{i+j}-S_i)^2 - \bkE (S_j^2) \, |
{\mathcal F}_0 \big ) \Vert_{1, \Phi, p} \, .
\end{eqnarray*}
Proceeding as in the proof of (4.6), p. 65 in Rio (2000), one can prove that,
for any $\sigma$-field ${\mathcal A}$ and any integrable random variable $X$,
 $\Vert \bkE ( X | {\mathcal A} ) \Vert_{1, \Phi, p} \leq \Vert  X  \Vert_{1, \Phi, p}$. Hence
\begin{eqnarray*}
\Vert \bkE \big( ( S_{i+j}-S_i)^2 - \bkE (S_j^2) \, | {\mathcal F}_0 \big) \Vert_{1, \Phi, p}
\leq
\Vert  \bkE (( S_{i+j}-S_i)^2 - \bkE (S_j^2) \, | {\mathcal F}_i \big ) \Vert_{1, \Phi, p} \, .
\end{eqnarray*}
By stationarity, it follows that $A_{i+j} \leq A_i + A_j$. Similarly $B_{i+j} \leq
B_i + B_j$. The proof of the equivalences then follows by using the
same arguments as in the proof of Lemma 2.7 in Peligrad and Utev
(2005).
\end{Remark}

\section{Rates of convergence for stationary sequences}
\label{Sectionstatseq}
\setcounter{equation}{0} In this section, we give estimates for the
ideal distances of order $r$ for stationary sequences which are not
necessarily adapted to ${\mathcal F}_i$.

\begin{Theorem} \label{propNAS} Let $(X_i)_{i \in \mathbb Z}$ be a stationary
sequence of centered random variables in ${\mathbb L}^p$ with $p \in
]2,3[$, and let $\sigma_n^2=n^{-1}{\mathbb E}(S_n^2)$. Assume that
\beq \label{Cond1cob}
 \sum_{n>0}{\mathbb E}(X_n | {\mathcal F}_0) \text{ and }
\sum_{n > 0} (X_{-n} - \bkE (X_{-n} | {\mathcal F}_0)) \text{ converge in ${\mathbb L}^p$} \, ,\eeq
and
 \beq
\label{Cond2cob} \sum_{n \geq 1}n^{-2+p/2} \Vert \,  n^{-1} \, \bkE
(S_n^2 | {\mathcal F}_0) - \sigma_n^2 \Vert_{p/2} < \infty \, .\eeq
Then the series $ \sum_{k \in {\mathbb Z}} \Cov(X_0 , X_k)$
converges to some nonnegative $\sigma^2$, and
\begin{enumerate}
\item $ \zeta_{r}( P_{n^{-1/2}S_n},G_{\sigma^2}
)=O(n^{1-p/2})$ for $r \in [p-2, 2]$, \item $\zeta_{r}(
P_{n^{-1/2}S_n},G_{\sigma_n^2} )=O(n^{1-p/2})$ for $r \in ]2, p]$.
\end{enumerate}
\end{Theorem}
\begin{Remark} According to the bound (\ref{cons1apMnSnp}), we infer
that,  under the assumptions of Theorem \ref{propNAS},  the
condition (\ref{Cond2cob}) is equivalent to \beq \label{Cond3cob}
\sum_{n \geq 1}n^{-2+p/2} \Vert \,  n^{-1} \,\bkE (S_n^2 | {\mathcal
F}_0) - \sigma^2 \Vert_{p/2} < \infty  \, .\eeq The same remark
applies to the next theorem with $p=3$.
\begin{Remark}
The result of item 1 is valid with $\sigma_n$ instead of $\sigma$.
On the contrary, the result of item 2 is no longer true if
$\sigma_n$ is replaced by $\sigma$, because for $r \in ]2, 3]$, a
necessary condition for  $\zeta_r(\mu,\nu)$ to be finite is that the
two first moments of $\nu$ and $\mu$ are equal. Note that  under the
assumptions of Theorem \ref{propNAS}, both $W_r(
P_{n^{-1/2}S_n},G_{\sigma^2})$ and  $W_r (
P_{n^{-1/2}S_n},G_{\sigma_n^2})$ are of the order of $n^{-(p-2)/ 2
\max (1, r)}$. Indeed, in the case where $r \in ]2, p]$, one has
that
$$
W_r ( P_{n^{-1/2}S_n},G_{\sigma^2})\leq W_r (
P_{n^{-1/2}S_n},G_{\sigma_n^2})+W_r (G_{\sigma_n^2},G_{\sigma^2})\,
,
$$
and the second term is of order $|\sigma-\sigma_n|=O(n^{-1/2})$.
\end{Remark}

\end{Remark}
In the case where $p=3$, the condition (\ref{Cond1cob}) has to be
strengthened.
\begin{Theorem} \label{propNASp3} Let $(X_i)_{i \in \mathbb Z}$ be a stationary
sequence of centered random variables in ${\mathbb L}^3$, and let
$\sigma_n^2=n^{-1}{\mathbb E}(S_n^2)$. Assume that \beq
\label{Condcobp3adap} \sum_{n \geq 1} \frac 1n \Big \Vert \sum_{k
\geq n } \bkE (X_k |{\mathcal F_0}) \Big \Vert_3 < \infty \quad
\text{and} \quad \sum_{n \geq 1} \frac 1n \Big \Vert \sum_{k \geq n
} ( X_{-k} - \bkE (X_{-k} |{\mathcal F_0}) ) \Big \Vert_3 < \infty
\, .\eeq Assume in addition that \beq \label{Cond2cobp3} \sum_{n
\geq 1}n^{-1/2}\Vert \,  n^{-1} \,\bkE (S_n^2 | {\mathcal F}_0) -
\sigma_n^2 \Vert_{3/2} < \infty \, . \eeq Then the series $ \sum_{k
\in {\mathbb Z}} \Cov(X_0 , X_k)$ converges to some nonnegative
$\sigma^2$ and
\begin{enumerate}
\item $\zeta_{1}( P_{n^{-1/2}S_n},G_{\sigma^2} )=O(n^{-1/2}\log n)$,
\item $\zeta_{r}( P_{n^{-1/2}S_n},G_{\sigma^2})=O(n^{-1/2})$ for $r \in ]1, 2]$,
\item $\zeta_{r}( P_{n^{-1/2}S_n},G_{\sigma_n^2} )=O(n^{-1/2})$ for $r \in ]2, 3]$.
\end{enumerate}
\end{Theorem}

\section{Applications} \label{SectionAppli}
\setcounter{equation}{0}
\subsection{Martingale differences sequences and  functions of Markov
chains}\label{MC}

Recall that the strong mixing coefficient of Rosenblatt (1956)
between two $\sigma$-algebras ${\mathcal A}$ and ${\mathcal B}$ is
defined by $
  \alpha ({\mathcal A}, {\mathcal B})=\sup \{ |{\mathbb P}(A \cap B)-{\mathbb P}(A){\mathbb
  P}(B)| \, : \, (A, B)  \in {\mathcal A} \times
  {\mathcal B} \, \} $.
 For a strictly stationary sequence $(X_i)_{i\in {\mathbb Z}}$, let
 ${\mathcal F}_i =\sigma (X_k, k \leq i)$. Define the mixing
 coefficients $\alpha_1(n)$ of the sequence $(X_i)_{i\in {\mathbb
 Z}}$ by
$$\alpha_1(n)=  \alpha ({\mathcal F}_0, \sigma(X_n))\,.$$ Let $Q$ be the quantile function of $|X_0|$, that is the
cadlag inverse  of the tail function  $x \rightarrow {\mathbb
P}(|X_0|>x)$. According to the results of Section 2, the following
proposition holds.

\begin{Proposition} \label{propalpha}
Let $(X_i)_{i \in \mathbb Z}$ be a stationary martingale difference
sequence. Assume moreover that the series \beq \label{condalpha1}
\sum_{k\geq 1}\frac{1}{k^{2-p/2}} \int_0^{\alpha_{1}(k)} ( 1 \vee
\log (1/u))^{(p-2)/2} Q^2(u) du \ \text{ and }\ \sum_{k\geq
1}\frac{1}{k^{2/p}} \Big (\int_0^{\alpha_{1}(k)} Q^p(u) du
\Big)^{2/p} \eeq are convergent.Then the conclusions of Theorem
\ref{Thmart} hold.
\end{Proposition}

\begin{Remark} From Theorem 2.1(b) in Dedecker and Rio (2007), a sufficient condition to get $W_1 (P_{n^{-1/2}S_n},
G_{\sigma^2})= O( n^{-1/2} \log n)$ is $$ \sum_{k \geq 0}
\int_0^{\alpha_1(n)} Q^3(u) du < \infty \, .$$ This condition is
always strictly stronger than the condition (\ref{condalpha1}) when
$p=3$. \end{Remark} We now give an example. Consider the homogeneous
Markov chain $(Y_i)_{i \in {\mathbb Z}}$ with state space ${\mathbb
Z}$ described at page 320 in Davydov (1973). The transition
probabilities are given by $p_{n,n+1}=p_{-n,-n-1}=a_n$ for $n \geq
0$, $p_{n,0}=p_{-n,0}=1-a_n$ for $n > 0$,  $p_{0,0}= 0$, $a_0 =1/2$
and $1/2 \leq  a_n < 1$ for $n \geq 1$. This chain is irreducible and
aperiodic.  It is Harris positively recurrent as soon as $\sum_{n
\geq 2}\Pi_{k=1}^{n-1}a_k < \infty$. In that case the stationary
chain is strongly mixing in the sense of Rosenblatt (1956).

Denote by  $K$ the Markov kernel of the chain  $(Y_i)_{i \in
{\mathbb Z}}$. The functions $f$ such that $K(f)=0$ almost
everywhere are obtained by linear combinations of the two functions
$f_1$ and $f_2$ given by $f_1(1)=1$, $f_1(-1)=-1$  and $f_1(n)=
f_1(-n)= 0$ if $n \neq 1$, and $f_2(0)=1$, $f_2(1)=f_2(-1)=0$ and
$f_2(n+1)=f_2(-n-1) = 1 - a_n^{-1}$ if $n
>0$. Hence the functions $f$ such that $K(f)=0$ are  bounded.

If  $(X_i)_{i \in \mathbb Z}$ is defined by $X_i = f(Y_i)$, with
$K(f)=0$, then Proposition  \ref{propalpha} applies if \beq
\label{condalphaharm} \alpha_1 (n) = O ( n^{1 -p/2} (\log n)^{-p/2 -
\epsilon} ) \, \text{ for some $\epsilon > 0$,} \eeq which holds as
soon as $P_0 (\tau = n ) =  O ( n^{-1 -p/2} (\log n)^{-p/2
-\epsilon} )$, where  $P_0$ is the probability of the chain starting
from $0$, and $\tau =\inf\{n>0, X_n=0\}$. Now $P_0 (\tau = n )=(1-
a_n)\Pi_{i=1}^{n-1}a_i$ for $n\geq 2$. Consequently, if
$$
a_i = 1 - \frac{p}{2i} \Big ( 1 + \frac{1 + \epsilon}{\log i} \Big )
\text{ for $i$ large enough}\, ,
$$
the condition (\ref{condalphaharm}) is satisfied and the conclusion
of Theorem \ref{Thmart} holds.

\begin{Remark} If $f$ is bounded and $K(f) \neq 0$,
the central limit theorem may fail to hold for $S_n = \sum_{i=1}^n
(f(Y_i) - \bkE (f(Y_i)))$. We refer to the Example 2, page 321,
given  by Davydov (1973), where $S_n$ properly normalized
converges to a stable law with exponent strictly less than 2.
\end{Remark}
\medskip

\noindent{\bf Proof of Proposition \ref{propalpha}.} Let
${B}^p({\mathcal F}_0)$ be the set of ${\mathcal F}_0$-measurable
random variables such that $\|Z\|_p  \leq 1$. We first notice that
$$\Vert \bkE (X^2_k |{\mathcal F}_0) - \sigma^2 \Vert_{p/2} = \sup_{Z \in {B}^{p/(p-2)}({\mathcal F}_0)} {\rm Cov}
(Z, X^2_k) \, .
$$
 Applying  Rio's covariance inequality
(1993), we get that
$$
\Vert \bkE (X^2_k |{\mathcal F}_0) - \sigma^2 \Vert_{p/2}  \leq 2
\Big ( \int_0^{\alpha_{1}(k)} Q^p(u) du \Big )^{2/p} \, ,
$$
which shows that the convergence of the second series in
(\ref{condalpha1})  implies (\ref{C2}). Now, from Fr\'echet  (1957),
we have that
$$\Vert \bkE (X^2_k |{\mathcal F}_0) - \sigma^2 \Vert_{1, \Phi, p} = \sup \big \{
\bkE ((1 \vee |Z|^{p-2}) | \, \bkE (X^2_k |{\mathcal F}_0)- \sigma^2
| \,  ) , Z \text{
 ${\mathcal F}_0$-measurable,
$Z \sim {\cal N} (0,1)$}\big \}\, .
$$
Hence, setting $\varepsilon_k = \text{sign} (\bkE (X^2_k |{\mathcal
F}_0)- \sigma^2)$,
$$\Vert \bkE (X^2_k |{\mathcal F}_0) - \sigma^2 \Vert_{1, \Phi, p} = \sup \big \{ {\rm Cov}
(\varepsilon_k (1 \vee |Z|^{p-2}), X^2_k), Z \text{
 ${\mathcal F}_0$-measurable,
$Z \sim {\cal N} (0,1)$}\big \}\, .
$$
Applying again  Rio's covariance inequality (1993), we get that
$$ \Vert \bkE (X^2_k |{\mathcal F}_0) - \sigma^2 \Vert_{1, \Phi, p} \leq
C \Big ( \int_0^{\alpha_{1}(k)} (1 \vee \log (u^{-1}))^{(p-2)/2}
Q^2(u) du \Big) \, ,
$$
which shows that the convergence of the first series in
(\ref{condalpha1}) implies (\ref{C1}).

\subsection{Linear processes and functions of linear processes}

\begin{Theorem} \label{propLP} Let $(a_i)_{i \in {\mathbb Z}}$ be a sequence of real numbers in  $\ell^2$
such that $\sum_{i \in {\mathbb Z}} a_i$ converges to some real $A$.
Let $(\varepsilon_i)_{i \in \mathbb Z}$ be a stationary sequence of
martingale differences in ${\mathbb L}^{p}$ for $p \in ]2,3]$. Let
$X_k=\sum_{j \in {\mathbb Z}}a_j \varepsilon_{k-j}$, and
$\sigma_n^2=n^{-1}{\mathbb E}(S_n^2)$. Let $b_0=a_0-A$ and $b_j=a_j$
for $j\neq 0$. Let $ A_n=\sum_{j \in {\mathbb Z}} ( \sum_{k=1}^n
b_{k-j})^2$. If $A_n=o(n)$, then $\sigma_n^2$ converges to
$\sigma^2=A^2 {\mathbb E}(\varepsilon_0^2)$. If moreover
\beq
\label{Condinnovations} \sum_{n =
1}^\infty\frac{1}{n^{2-p/2}}\Big \Vert \bkE \Big(
\frac{1}{n}\Big(\sum_{j=1}^n \varepsilon_j \Big)^2 \Big| {\mathcal
F}_0 \Big ) - {\mathbb E}(\varepsilon_0^2) \Big \Vert_{p/2} < \infty
\, ,
\eeq
then we have
\begin{enumerate}
\item  If $A_n= O(1)$, then $ \zeta_{1}( P_{n^{-1/2}S_n},G_{\sigma^2}
)=O(n^{-1/2}\log(n))$, for $p=3$,
\item If $A_n=O( n^{(r+2-p)/r})$, then
$ \zeta_{r}( P_{n^{-1/2}S_n},G_{\sigma^2} )=O(n^{1-p/2})$, for $r
\in [p-2, 1]$ and $p\neq 3$,
\item If $A_n=O( n^{3-p})$, then
$ \zeta_{r}( P_{n^{-1/2}S_n},G_{\sigma^2} )=O(n^{1-p/2})$, for $r
\in ]1, 2]$,
 \item
 If $A_n=O( n^{3-p})$, then
$\zeta_{r}( P_{n^{-1/2}S_n},G_{\sigma_{n}^2} )=O(n^{1-p/2})$, for $r
\in ]2, p]$.
\end{enumerate}
\end{Theorem}

\begin{Remark}\label{Rempl}
If the condition given by Heyde (1975) holds, that is
\begin{equation}\label{heyde}
\sum_{n=1}^{\infty }\Big( \sum_{k\geq n}a_{k}\Big) ^{2}<\infty \quad
\text{and} \quad \sum_{n=1}^{\infty }\Big( \sum_{k\leq -
n}a_{k}\Big) ^{2}<\infty \, ,
\end{equation}
 then $A_n=O(1)$, so that it satisfies all the conditions of items
1-4. On the other and, one has the bound
\begin{equation}\label{alternatheyde}
 A_n \leq 4 B_n, \quad \text{where} \quad  B_n=\sum_{k=1}^n \Big (  \Big(\sum_{j\geq k} |a_j|\Big)^2 + \Big(\sum_{j \leq -k} |a_j|\Big)^2\Big)\, .
\end{equation}

\end{Remark}

\noindent{\bf Proof of Theorem \ref{propLP}.} We start with the
following decomposition:
\beq \label{declinear} S_n = A \sum_{j=1}^n
\varepsilon_j + \sum_{j= - \infty}^{\infty} \Big ( \sum_{k=1}^n
b_{k-j}\Big ) \varepsilon_j \, .
\eeq
Let $R_n = \sum_{j= -\infty}^{\infty} ( \sum_{k=1}^n b_{k-j} ) \varepsilon_j $. Since
$\|R_n\|_2^2=A_n \|\varepsilon_0\|_2^2$ and since $|\sigma_n-\sigma|
\leq n^{-1/2}\|R_n\|_2$, the fact that $A_n=o(n)$ implies that
$\sigma_n$ converges to $\sigma$. We now give an upper bound
for $\|R_n\|_p$. From Burkholder's inequality, there exists a
constant $C$ such that
\begin{equation}\label{burk}
\Vert R_n \Vert_p  \leq  C \Big \{ \Big\Vert \sum_{j= -
\infty}^{\infty}  \Big ( \sum_{k=1}^n b_{k-j}\Big )^2
\varepsilon_j^2 \Big\Vert_{p/2} \Big \}^{1/2}  \leq  C \Vert
\varepsilon_0\Vert_{p} \sqrt{A_n}.
\end{equation}The result follows by  applying Theorem \ref{Thmart} to the martingale
$A\sum_{k=1}^n \varepsilon_k$ (this is possible because of
(\ref{Condinnovations})), and by using Lemma \ref{compSnMn}  with
the upper bound (\ref{burk}).  To prove Remark \ref{Rempl}, note
first that $$ A_n= \sum_{j=1}^n \Big ( \sum_{l=-\infty}^{-j} a_l +
\sum_{l=n+1-j}^\infty a_l\Big )^2 + \sum_{i=1}^\infty \Big(
\sum_{l=i}^{n+i-1} a_l \Big)^2 + \sum_{i=1}^\infty \Big (
\sum_{l=-i-n+1}^{-i}  a_l \Big)^2 \, .
$$
It follows easily that $A_n=O(1)$ under (\ref{heyde}). To prove the
bound (\ref{alternatheyde}), note first that $$ A_n \leq 3 B_n +
\sum_{i=n+1}^\infty \Big( \sum_{l=i}^{n+i-1} |a_l|
\Big)^2+\sum_{i=n+1}^\infty \Big ( \sum_{l=-i-n+1}^{-i}  |a_l|
\Big)^2\, .$$ Let $T_i=\sum_{l=i}^\infty |a_l|$ and $Q_i=
\sum_{l=-\infty}^{-i} |a_l|$. We have that
\begin{eqnarray*}
\sum_{i=n+1}^\infty \Big( \sum_{l=i}^{n+i-1} |a_l| \Big)^2&\leq&
T_{n+1}\sum_{i=n+1}^\infty (T_i-T_{n+i}) \leq n T_{n+1}^2 \\
\sum_{i=n+1}^\infty \Big( \sum_{l=-i-n+1}^{-i}  |a_l|\Big)^2&\leq&
Q_{n+1}\sum_{i=n+1}^\infty (Q_i-Q_{n+i}) \leq n Q_{n+1}^2.
\end{eqnarray*}
Since $n(T_{n+1}^2+ Q_{n+1}^2) \leq B_n$, (\ref{alternatheyde})
follows. $\square$

\medskip

In the next result, we shall focus on functions of real-valued
linear processes
\begin{equation}\label{def2suite}
X_k=h\Big(\sum_{i \in {\mathbb Z}} a_i \varepsilon_{k-i}\Big)-
{\bkE}\Big(h\Big(\sum_{i \in {\mathbb Z}} a_i
\varepsilon_{k-i}\Big)\Big) \, ,
\end{equation}
where $(\varepsilon_i)_{i \in {\mathbb Z}}$ is  a sequence of iid
random variables. Denote by $w_h(.,M)$ the modulus of continuity of
the function $h$ on the interval $[-M,M]$, that is
$$
  w_h(t,M)=\sup  \{ |h(x)-h(y)|, |x-y| \leq t, |x|\leq M, |y|
  \leq M  \} \, .
$$

\begin{Theorem}\label{Thlin} Let
$(a_i)_{i \in {{\mathbb Z}}}$ be a sequence of real numbers in
$\ell^2$ and $(\varepsilon_i)_{i \in \mathbb Z}$ be a sequence of
iid random variables in ${\mathbb L}^{2}$. Let $X_k$  be defined as
in (\ref{def2suite}) and $\sigma_n^2=n^{-1}{\mathbb E}(S_n^2)$.
Assume that $h$ is $\gamma$-H\"older on any compact set, with
$w_h(t, M) \leq C t^\gamma M^\alpha$, for some $C>0$, $\gamma \in
]0,1]$ and $\alpha \geq 0$. If for some $p \in ]2, 3]$,
\begin{equation}\label{gammaholder}
{\mathbb E}(|\varepsilon_0|^{2 \vee (\alpha+\gamma)p})<\infty \quad
\text{and} \quad \sum_{i \geq 1} i^{p/2-1}\Big( \sum_{|j|\geq i}
a_j^2 \Big)^{\gamma/2} < \infty,
\end{equation}
then the series $\sum_{k \in {\mathbb Z}}\Cov(X_0,X_k)$ converges to
some nonnegative $\sigma^2$, and
\begin{enumerate}
\item $ \zeta_{1}( P_{n^{-1/2}S_n},G_{\sigma^2}
)=O(n^{-1/2}\log n)$, for $p=3$,
\item $ \zeta_{r}( P_{n^{-1/2}S_n},G_{\sigma^2}
)=O(n^{1-p/2})$ for $r \in [p-2, 2]$ and $(r,p)\neq (1,3)$,
\item $\zeta_{r}( P_{n^{-1/2}S_n},G_{\sigma_n^2} )=O(n^{1-p/2})$ for $r \in ]2, p]$.
\end{enumerate}
\end{Theorem}

\noindent {\bf Proof of Theorem \ref{Thlin}.} Theorem \ref{Thlin} is
a consequence of the following proposition:
\begin{Proposition} \label{modulo}
Let $(a_i)_{i \in {\mathbb Z}}$, $(\varepsilon_i)_{i \in {\mathbb
Z}}$ and $(X_i)_{i \in {\mathbb Z}}$ be as in Theorem \ref{Thlin}.
Let $(\varepsilon'_i)_{i \in {\mathbb Z}}$ be an independent copy of
$(\varepsilon_i)_{i \in {\mathbb Z}}$. Let $V_0= \sum_{i \in
{\mathbb Z}} a_i \varepsilon_{-i}$ and $$M_{1,i}= |V_0| \vee
\Big|\sum_{j<i}a_j \varepsilon_{-j}+\sum_{j \geq i} a_j
\varepsilon'_{-j}\Big|\quad \text{and} \quad  M_{2, i}=|V_0| \vee
\Big|\sum_{j< i}a_j \varepsilon'_{-j}+\sum_{j \geq i} a_{j}
\varepsilon_{-j}\Big|.$$  If for some $p \in ]2,3]$, \beq
\label{mod}
 \sum_{i \geq 1} i^{p/2-1} \Big \| w_h\Big(\Big|\sum_{j \geq i} a_j
 \varepsilon_{-j}\Big |, M_{1,i}\Big) \Big \|_p < \infty \quad \text{and}
 \quad
 \sum_{i \geq 1} i^{p/2-1} \Big \| w_h\Big(\Big|\sum_{j <-i} a_j
 \varepsilon_{-j}\Big|, M_{2,-i}\Big) \Big \|_p < \infty,
\eeq then the conclusions
 of Theorem \ref{Thlin} hold.
\end{Proposition}
\noindent To prove Theorem \ref{Thlin}, it remains to check
(\ref{mod}). We only check the first condition. Since $w_h(t, M)
\leq C t^\gamma M^\alpha$ and the random variables $\varepsilon_i$
are iid, we have
\begin{eqnarray*}
\Big \| w_h\Big(\Big |\sum_{j \geq i} a_j
 \varepsilon_{-j}\Big|, M_{1,i}\Big) \Big \|_p &\leq & C\Big \|\Big |\sum_{j
 \geq i} a_j \varepsilon_{-j} \Big |^\gamma |V_0|^\alpha \Big \|_p + C\Big \|\Big |\sum_{j
 \geq i} a_j \varepsilon_{-j} \Big |^\gamma \Big \|_p \| |V_0|^\alpha
 \|_p\, ,
 \end{eqnarray*}
 so that
 \begin{multline*}
\Big \| w_h\Big(\Big |\sum_{j \geq i} a_j
 \varepsilon_{-j}\Big|, M_{1,i}\Big) \Big \|_p \\\leq
 C \Big ( 2^\alpha \Big\|\Big |\sum_{j
 \geq i} a_j \varepsilon_{-j} \Big |^{\alpha+ \gamma}\Big \|_p +
 \Big\|\Big |\sum_{j
 \geq i} a_j \varepsilon_{-j} \Big |^{\gamma}\Big \|_p\Big(\|
 |V_0|^\alpha
 \|_p+2^\alpha \Big\|\Big |\sum_{j
 < i} a_j \varepsilon_{-j} \Big |^{\alpha}\Big \|_p \Big)\Big).
\end{multline*}
From Burkholder's inequality, for any $\beta>0$, $$ \Big\|\Big
|\sum_{j
 \geq i} a_j \varepsilon_{-j} \Big |^{\beta}\Big \|_p=\Big\|\sum_{j
 \geq i} a_j \varepsilon_{-j}\Big \|_{\beta p}^\beta \leq K
 \Big (\sum_{j \geq i}
 a_j^2\Big)^{\beta/2}\|\varepsilon_0\|_{2\vee\beta p}^{\beta} \,
 .
$$
Applying this inequality with $\beta=\gamma$ or $\beta= \alpha+
\gamma$, we infer that the first part of (\ref{mod}) holds under
(\ref{gammaholder}). The second part can be handled in the same way.
$\square$

\medskip

\noindent {\bf Proof of Proposition \ref{modulo}.} Let ${\mathcal
F}_i=\sigma (\varepsilon_k, k \leq i)$. We shall first prove that
the condition (\ref{Cond2cob}) of Theorem \ref{propNAS} holds. We
write
\begin{eqnarray*}
\Vert \bkE (S_n^2 | {\mathcal F}_0)\!\! &-& \!\!\bkE (S_n^2)
\Vert_{p/2} \leq
 2 \sum_{i=1}^n\sum_{k=0}^{n-i}\Vert \bkE (X_iX_{k+i} |
{\mathcal F}_0) -  \bkE (X_iX_{k+i}) \Vert_{p/2} \\
& \leq & 4 \sum_{i=1}^{n}\sum_{k=i}^{n}\Vert \bkE (X_iX_{k+i}  |
{\mathcal F}_0) \Vert_{p/2}
 + 2 \sum_{i=1}^n\sum_{k=1}^{i} \Vert
\bkE (X_iX_{k+i} | {\mathcal F}_0) - \bkE (X_iX_{k+i}) \Vert_{p/2}\,
.
\end{eqnarray*}
We first control the second term. Let $\varepsilon^{\prime}$ be an
independent copy of $\varepsilon$, and denote by
$\bkE_{\varepsilon}(\cdot)$ the conditional expectation with respect
to $\varepsilon$. Define
\begin{eqnarray*}
Y_i = \sum_{j <i} a_j \varepsilon_{i-j}\;,\; Y_i^\prime = \sum_{j
<i} a_j \varepsilon^{\prime}_{i-j}\;, Z_i = \sum_{j \geq i} a_j
\varepsilon_{i-j}\;,\; Z_i^\prime = \sum_{j \geq i} a_j
\varepsilon^{\prime}_{i-j}\;
\end{eqnarray*}
and $m_{1,i}=|Y'_i+Z_i|\vee |Y'_i+Z'_i|$. Taking
$\mathcal{F}_{\ell}=\sigma ( \varepsilon_i, i \leq \ell)$, and
setting $h_0 = h - \bkE ( h ( \sum_{i \in {\mathbb Z}} a_i
\varepsilon_i))$, we have
\begin{multline*}
  \| \bkE ( X_iX_{k+i} |\mathcal{F}_{0})-\bkE (X_iX_{k+i})  \|_{p/2} \\
 = \Big \|{\bkE}_{\varepsilon}\Big(h_0(Y_i^\prime + Z_i) h_0%
(Y_{k+i}^\prime + Z_{k+i} )\Big) -
{\bkE}_{\varepsilon}\Big(h_0(Y_i^\prime + Z_i^\prime) h_0%
(Y_{k+i}^\prime + Z_{k+i}^\prime )\Big) \Big \|_{p/2} \, .
\end{multline*}
Hence,
\begin{eqnarray*}
  \| \bkE ( X_iX_{k+i} |\mathcal{F}_{0})-\bkE (X_iX_{k+i})  \|_{p/2}
  &  \leq & \Vert h_0 (Y_{k+i}^\prime + Z_{k+i} ) \Vert_p \Big \Vert
w_h\Big(\Big|\sum_{j \geq i} a_j (\varepsilon_{i-j}-
\varepsilon_{i-j}^\prime\Big)\Big|, m_{1,i} \Big)\Big \Vert_p \\
&+&\!\!\! \Vert h_0 (Y_{i}^\prime + Z_{i}^\prime )\Vert_p \Big\Vert
w_h\Big(\Big|\sum_{j \geq k+i} a_j (\varepsilon_{k+i-j} -
\varepsilon_{k+i-j}^\prime\Big)\Big|, m_{1, k+i} \Big)\Big\Vert_p \,
.
\end{eqnarray*}
By subadditivity,
\begin{eqnarray*}
\Big\|w_h\Big(\Big|\sum_{j \geq i} a_j (\varepsilon_{i-j}-
\varepsilon_{i-j}^\prime)\Big|, m_{1,i} \Big)\Big \|_p&\leq &\Big
\|w_h \Big(\Big |\sum_{j \geq i} a_j \varepsilon_{i-j}\Big |,
m_{1,i} \Big )\Big \|_p+\Big \|w_h\Big (\Big|\sum_{j \geq i} a_j
\varepsilon_{i-j}^\prime\Big |, m_{1,i} \Big)\Big\|_p\\
&\leq & 2 \Big \| w_h\Big(\Big|\sum_{j \geq i} a_j
 \varepsilon_{-j}\Big |, M_{1,i}\Big) \Big \|_p \, .
\end{eqnarray*}
In the same way
$$
\Big \Vert w_h\Big(\Big|\sum_{j \geq k+i} a_j (\varepsilon_{k+i-j} -
\varepsilon_{k+i-j}^\prime)\Big|, m_{1, k+i}\Big)\Big \Vert_p \leq 2
\Big \| w_h\Big(\Big|\sum_{j \geq k+i} a_j
 \varepsilon_{-j}\Big |, M_{1,k+i}\Big) \Big \|_p \, .
$$
Consequently $$ \sum_{n \geq 1} \frac{1}{n^{3 - p/2}}
\sum_{i=1}^n\sum_{k=1}^{i} \Vert \bkE (X_iX_{k+i} | {\mathcal F}_0)
- \bkE (X_iX_{k+i}) \Vert_{p/2} < \infty $$ provided that the first
condition in (\ref{mod}) holds.

We turn now to the control of $\sum_{i=1}^{n}\sum_{k=i}^{n}\Vert
\bkE (X_iX_{k+i}  | {\mathcal F}_0) \Vert_{p/2}$. We first write
that
\begin{eqnarray*}
 \| \bkE ( X_iX_{k+i} |\mathcal{F}_{0}) \|%
_{p/2} & = &  \| \bkE \big ( (X_i - \bkE ( X_i |\mathcal{F}_{i+[k/2]}))X_{k+i} |\mathcal{F}_{0} \big ) \|%
_{p/2} +  \| \bkE \big ( \bkE ( X_i |\mathcal{F}_{i+[k/2]})X_{k+i} |\mathcal{F}_{0} \big )\|%
_{p/2} \\
& = & \Vert X_0 \Vert_p  \| X_i - \bkE ( X_i |\mathcal{F}_{i+[k/2]})
\|_p
+ \Vert X_0 \Vert_p  \|\bkE (X_{k+i} |\mathcal{F}_{i+[k/2]})\|%
_{p} \, .
\end{eqnarray*}
Let $b(k)=k-[k/2]$. Since $ \|\bkE (X_{k+i} |\mathcal{F}_{i+[k/2]})
\| _{p}  =  \|\bkE (X_{b(k)} |\mathcal{F}_0) \| _{p}$, we have that
\begin{multline*}  \label{bof}
 \|\bkE (X_{k+i} |\mathcal{F}_{i+[k/2]}) \| _{p}
\\
=\Big \|{\bkE}_{\varepsilon}\Big(h\Big(\sum_{j <b(k)}
a_j \varepsilon^{\prime}_{b(k)-j} + \sum_{j \geq b(k)} a_j \varepsilon_{b(k)-j}\Big) %
- h\Big(\sum_{j <b(k)} a_j \varepsilon^{\prime}_{b(k)-j} + \sum_{j
\geq b(k)} a_j \varepsilon^{\prime}_{b(k)-j}\Big)  \Big) \Big \|%
_{p}\, .
\end{multline*}
Using the same arguments  as before, we get that
$$
 \|\bkE (X_{k+i} |\mathcal{F}_{i+[k/2]}) \| _{p} \leq 2 \Big \|
w_h\Big(\Big|\sum_{j \geq b(k)} a_j
 \varepsilon_{-j}\Big |, M_{1,b(k)}\Big) \Big \|_p\, .
$$
In the same way,
\begin{multline*}
\Big \| X_i - \bkE ( X_i |\mathcal{F}_{i+[k/2]})\Big \|_p  \\ =
 \Big \|{\bkE}_{\varepsilon}\Big(h\Big(\sum_{j <-[k/2]}
a_j \varepsilon_{i-j} + \sum_{j \geq -[k/2]} a_j
\varepsilon_{i-j}\Big) - h\Big(\sum_{j <-[k/2]} a_j
\varepsilon^{\prime}_{i-j} + \sum_{j
\geq -[k/2]} a_j \varepsilon_{i-j}\Big)  \Big) \Big \|%
_{p}  \, ,
\end{multline*}
so that
$$
\Big \| X_i - \bkE ( X_i |\mathcal{F}_{i+[k/2]})\Big \|_p \leq
2\Big \| w_h\Big(\Big|\sum_{j <-[k/2]} a_j
 \varepsilon_{-j}\Big |, M_{2,-[k/2]}\Big) \Big \|_p \, .
 $$
  Consequently
$$ \sum_{n \geq 1} \frac{1}{n^{3 - p/2}}
\sum_{i=1}^{n}\sum_{k=i}^{n}\Vert \bkE (X_iX_{k+i}  | {\mathcal
F}_0) \Vert_{p/2} < \infty $$ provided that (\ref{mod}) holds. This
completes the proof of (\ref{Cond2cob}). Using the same arguments,
one can easily check that the condition (\ref{Cond1cob}) of Theorem
\ref{propNAS} (and also  the condition (\ref{Condcobp3adap}) of
Theorem  \ref{propNASp3} in the case $p=3$) holds under (\ref{mod}).
 $\square$

\subsection{Functions of $\phi$-dependent sequences}

In order to include examples of dynamical systems satisfying some
correlations inequalities, we introduce a weak version of the
uniform mixing coefficients (see Dedecker and Prieur (2007)).

\begin{Definition}
For any random variable $Y=(Y_1, \cdots, Y_k)$ with values in
${\mathbb R}^k$ define the function  $g_{x,j}(t)=\I_{t \leq
x}-{\mathbb P}(Y_j \leq x)$. For any $\sigma$-algebra ${\cal F}$,
let
$$
\phi({\cal F}, Y)= \sup_{(x_1, \ldots , x_k) \in {\mathbb R}^k}
\Big \|{\mathbb E}\Big (\prod_{j=1}^k g_{x_j,j}(Y_j)\Big |{\cal F}
\Big)- {\mathbb E}\Big (\prod_{j=1}^k
g_{x_j,j}(Y_j)\Big)\Big\|_{\infty}.
$$
For a sequence ${\bf Y}=(Y_i)_{i \in {\mathbb Z}}$, where $Y_i=Y_0
\circ T^i$ and $Y_0$ is a ${\cal F}_0$-measurable and real-valued
r.v., let
$$
\phi_{k, {\bf Y}}(n) = \max_{1 \leq l \leq k} \ \sup_{ i_l>\ldots
> i_1 \geq n} \phi({\cal F}_0, (Y_{i_1}, \ldots, Y_{i_l})) .
$$
\end{Definition}
\begin{Definition} \label{defclosedenv}
For any $p \geq 1$, let ${\mathcal C} (p, M, P_X )$ be the closed
convex envelop of the set of  functions $f$ which are monotonous on
some open interval of ${\mathbb R}$ and null elsewhere, and such
that $\bkE(|f(X)|^p) < M$.
\end{Definition}

\begin{Proposition} \label{propphimixing} Let $p \in ]2,3]$ and $s \geq p$. Let  $X_i = f(Y_i) - \bkE ( f(Y_i))$, where
$Y_i=Y_0 \circ T^i$ and $f$ belongs to ${\mathcal C}(s,M, P_{Y_0})$.
Assume that \beq \label{condphi} \sum_{i\geq 1}i^{(p-4)/2 +
(s-2)/(s-1)} \phi_{2, {\bf Y}}(i)^{(s-2)/s} < \infty \, . \eeq Then
the conclusions
 of Theorem \ref{Thlin}
hold.
\end{Proposition}
\begin{Remark}
Notice that if $s=p=3$, the condition (\ref{condphi}) becomes $
\sum_{i\geq 1} \phi_{2, {\bf Y}}(i) ^{1/3} < \infty \, , $ and if $s
= \infty$, the condition (\ref{condphi}) becomes $ \sum_{i\geq
1}i^{(p-2)/2} \phi_{2, {\bf Y}}(i)  < \infty. $
\end{Remark}

\noindent {\bf Proof of Proposition \ref{propphimixing}.} Let
${B}^p({\mathcal F}_0)$ be the set of ${\mathcal F}_0$-measurable
random variables such that $\|Z\|_p  \leq 1$. We first notice that
$$ \|\bkE (X_k |{\mathcal F}_0)\Vert_{p} \leq \Vert \bkE (X_k
|{\mathcal F}_0) \Vert_{s}= \sup_{Z \in B^{s/(s-1)}({\mathcal F}_0)}
{\rm Cov} (Z, f(Y_k)) \, .
$$
According to Corollary \ref{propineq3} and since $ \phi(\sigma(Z),
Y_k) \leq \phi_{1, {\bf Y}}(k) \, , $ we get that \begin{equation}
\label{app1covphi}\Vert \bkE (X_k |{\mathcal F}_0) \Vert_{s} \leq 8
M^{1/s} (\phi_{1, {\bf Y}}(k) )^{(s-1)/s} \, .\end{equation} It
follows that the conditions (\ref{Cond1cob}) (for $p \in ]2,3[$) or
(\ref{Condcobp3adap}) (for $p=3$) are satisfied under
(\ref{condphi}).  The condition (\ref{Cond2cob}) follows from the
following lemma by taking $b=(4-p)/2$.
\begin{Lemma} \label{b} Let $X_i$ be as in Proposition \ref{propphimixing},
and let $b \in ]0, 1[$.   $$ \text{If} \quad \sum_{i\geq 1}i^{-b +
(s-2)/(s-1)} \phi_{2, {\bf Y}}(i)^{(s-2)/s} < \infty, \quad
\text{then} \quad \sum_{n >1} \frac{1}{n^{1+b}}\|\bkE (S_n^2 |
{\mathcal F}_0) - \bkE (S_n^2) \Vert_{p/2} < \infty \, .
$$
\end{Lemma}
\noindent {\bf Proof of Lemma \ref{b}.} Since,
$$
 \Vert \bkE (S_n^2 | {\mathcal F}_0) - \bkE (S_n^2) \Vert_{p/2}
 \leq 2 \sum_{i=1}^n\sum_{k=0}^{n-i}\Vert \bkE (X_iX_{k+i} |
{\mathcal F}_0) - \bkE (X_iX_{k+i}) \Vert_{p/2}, $$ we infer that
there exists $C>0$ such that
\begin{equation}\label{b1}
\sum_{n >1} \frac{1}{n^{1+b}}\|\bkE (S_n^2 | {\mathcal F}_0) - \bkE
(S_n^2) \Vert_{p/2} \leq C \sum_{i>0} \sum_{k\geq 0}
\frac{1}{(i+k)^b} \Vert \bkE (X_iX_{k+i} | {\mathcal F}_0) - \bkE
(X_iX_{k+i}) \Vert_{p/2}\, .
\end{equation}
We shall bound up $\Vert \bkE (X_iX_{k+i} | {\mathcal F}_0) - \bkE
(X_iX_{k+i}) \Vert_{p/2}$ in two ways. First, using the stationarity
and the upper bound (\ref{app1covphi}), we have that
\begin{equation*}
\Vert \bkE (X_iX_{k+i} | {\mathcal F}_0) - \bkE (X_iX_{k+i})
\Vert_{p/2} \leq 2 \|X_0{\mathbb E}(X_k|{\mathcal F}_0)\|_{p/2}\leq
16 \|X_0\|_p M^{1/s} (\phi_{1, {\bf Y}}(k) )^{(s-1)/s}\, .
\end{equation*}
Next, using again Corollary \ref{propineq3},
$$
\Vert \bkE (X_iX_{k+i} | {\mathcal F}_0) - \bkE (X_iX_{k+i})
\Vert_{p/2}\leq \sup_{Z \in B^{s/(s-2)}({\mathcal F}_0)} {\rm Cov}
(Z, X_iX_{k+i}) \leq 32 M^{2/s} (\phi_{2, {\bf Y}}(i) )^{(s-2)/s} \,
.
$$
From (\ref{b1}) and the above upper bounds, we infer that the
conclusion of Lemma (\ref{b}) holds provided that
$$ \sum_{i>0} \Big( \sum_{k=1}^{[i^{(s-2)/(s-1)}]} \frac{1}{(i+k)^b}
\Big)(\phi_{2, {\bf Y}}(i) )^{(s-2)/s} + \sum_{k\geq 0} \Big(
\sum_{i=1}^{[k^{(s-1)/(s-2)}]} \frac{1}{(i+k)^b} \Big)(\phi_{1, {\bf
Y}}(k) )^{(s-1)/s} < \infty \, .
$$
Here, note that
$$
\sum_{k=1}^{ [i^{(s-2)/(s-1)}]} \frac{1}{(i+k)^b}\leq
i^{-b+\frac{s-2}{s-1}}\quad \text{and} \quad
\sum_{i=1}^{[k^{(s-1)/(s-2)}]} \frac{1}{(i+k)^b}\leq
\sum_{m=1}^{[2k^{(s-1)/(s-2)}]} \frac{1}{m^b} \leq D
k^{(1-b)\frac{(s-1)}{(s-2)}}\, ,
$$
for some $D>0$. Since $\phi_{1, {\bf Y}}(k)\leq \phi_{2, {\bf
Y}}(k)$, the conclusion of lemma (\ref{b}) holds provided
$$
\sum_{i\geq 1}i^{-b + \frac{s-2}{s-1}} \phi_{2, {\bf
Y}}(i)^{\frac{s-2}{s}} < \infty \quad \text{and} \quad \sum_{k\geq
1}k^{(1-b)\frac{(s-1)}{(s-2)}} \phi_{2, {\bf Y}}(k)^{\frac{s-1}{s}}
< \infty \, .
$$
One can prove that the second series converges provided the first
one does.  $\square$

\subsubsection{Application to Expanding maps}
 Let $BV$ be the class of bounded variation functions from $[0,1]$
to ${\mathbb R}$. For any $h \in BV$, denote by $\|dh \|$ the
variation norm of the measure $dh$.

Let $T$ be a map from $[0, 1]$ to $[0, 1]$ preserving a probability
$\mu$ on $[0, 1]$,  and let
\begin{eqnarray*}
  S_n (f)=
\sum_{k=1}^{n} (f \circ T^{k}- \mu(f)) \, .
\end{eqnarray*}
Define  the Perron-Frobenius operator $K$ from ${\mathbb L}^2([0,
1], \mu)$ to ${\mathbb L}^2([0, 1], \mu)$ $via$ the equality
\begin{equation}\label{pfeq}
   \int_0^1 (Kh)(x) f(x) \mu (dx) = \int_0^1 h(x) (f \circ T)(x)
   \mu(dx)
   \, .
\end{equation}
A Markov Kernel $K$ is said to be $BV$-contracting if there exist
$C>0$ and $\rho \in [0, 1[$ such that
\begin{equation}\label{expan}
     \|d K^n(h)\| \leq C \rho^n \|dh \| \, .
\end{equation}
The map $T$ is said to be $BV$-contracting if its Perron-Frobenius
operator is $BV$-contracting.

Let us present a large class of $BV$-contracting maps. We shall say
that  $T$ is uniformly expanding if it  belongs to the class
${\mathcal C}$ defined in Broise (1996), Section 2.1 page 11. Recall
that if $T$ is uniformly expanding, then there exists a probability
measure $\mu$ on $[0, 1]$, whose density $f_\mu$ with respect to the
Lebesgue measure is a bounded variation function, and such that
$\mu$ is invariant by $T$. Consider now the more restrictive
conditions:
\begin{enumerate}
\item[(a)] $T$ is uniformly expanding.
\item[(b)] The invariant measure $\mu$ is unique and $(T, \mu)$ is
mixing in the ergodic-theoretic sense.
\item[(c)] $\displaystyle \frac{1}{f_\mu}{\bf 1}_{f_\mu>0}$ is a
bounded variation function.
\end{enumerate}
Starting from Proposition 4.11 in Broise (1996), one can prove that
if $T$ satisfies the assumptions (a), (b) and (c) above, then it is
$BV$ contracting (see for instance Dedecker and Prieur  (2007),
Section 6.3). Some well known examples of maps satisfying the
conditions (a), (b) and (c) are:
\begin{enumerate}
\item $T(x)= \beta x -[\beta x]$ for $\beta> 1$. These maps are called
$\beta$-transformations.
\item  $I$ is the finite union of disjoint intervals $(I_k)_{1 \leq k
\leq n}$, and $T(x)=a_kx +b_k$ on $I_k$, with $|a_k|>1$.
\item $T(x)=a(x^{-1}-1)-[a(x^{-1}-1)]$ for some $a>0$. For $a=1$, this
transformation is known as the Gauss map.
\end{enumerate}

\begin{Proposition} \label{appexpmap}
Let $\sigma_n^2=n^{-1}{\mathbb E}(S_n^2(f))$. If $T$ is
$BV$-contracting, and if $f$ belongs to $ {\mathcal C}(p,M, \mu )$
with $p \in ]2,3]$, then the series $ \mu((f -\mu(f))^2) + 2 \sum_{n
> 0} \mu(f \circ T^n \cdot (f-\mu(f)))$ converges to some
nonnegative $\sigma^2$, and
\begin{enumerate}
\item $ \zeta_{1}( P_{n^{-1/2}S_n(f)},G_{\sigma^2}
)=O(n^{-1/2}\log n )$, for $p=3$,
\item $ \zeta_{r}( P_{n^{-1/2}S_n(f)},G_{\sigma^2}
)=O(n^{1-p/2})$ for $r \in [p-2, 2]$ and  $(r,p) \neq (1,3)$,
 \item
$\zeta_{r}( P_{n^{-1/2}S_n(f)},G_{\sigma_{n}^2} )=O(n^{1-p/2})$ for
$r \in ]2, p]$.
\end{enumerate}
\end{Proposition}

\noindent {\bf Proof of Proposition \ref{appexpmap}.}
 Let $(Y_i)_{i \geq 1}$ be the Markov chain with transition Kernel
$K$ and invariant measure $\mu$. Using the equation (\ref{pfeq}) it
is easy to see that $(Y_0, \ldots, Y_n)$ it is distributed as
$(T^{n+1}, \ldots, T)$. Consequently, to prove Proposition
\ref{appexpmap}, it suffices to prove that the sequence $X_i=
f(Y_i)-\mu(f)$ satisfies the condition (\ref{condphi}) of
Proposition \ref{propphimixing}.

According to Lemma 1 in Dedecker and Prieur (2007), the coefficients
$\phi_{2, {\bf Y}} (i)$ of the chain $(Y_i)_{i \geq 0}$ with respect
to ${\mathcal F}_i = \sigma ( Y_j , j \leq i)$ satisfy
$\phi_{2, {\bf Y}} (i) \leq C \rho^i$
for some $\rho \in ]0, 1[$ and some positive constant $C$.
It follows that (\ref{condphi}) is satisfied for $s = p$.

\section{Proofs of the main results} \label{SectionProofs}
\setcounter{equation}{0}

From now on, we denote by $C$ a numerical constant  which may vary
from line to line.
\begin{Notation}\label{seminorm}
For $l$ integer, $q$ in $]l,l+1]$ and $f$ $l$-times continuously
differentiable, we set
$$
|f|_{\Lambda_q}   =  \sup
\{ |x-y|^{l-q} |f^{(l)} (x) - f^{(l)} (y)| : (x,y) \in {\mathbb R} \times {\mathbb R} \} .
$$
\end{Notation}

\subsection{Proof of Theorem \ref{Thmart}}

We prove Theorem \ref{Thmart} in the case $\sigma = 1$. The general
case follows by dividing the random variables by $\sigma$.  Since
$\zeta_r(P_{aX}, P_{aY})= |a|^r \zeta_r(P_{X}, P_{Y})$, it is enough
to bound up $\zeta_r(P_{S_n}, G_n)$. We first give an upper bound
for $ \zeta_{p,N} := \zeta_p ( P_{S_{2^N}} , G_{2^N} )$.
\par\msk\no
\begin{Proposition} \label{propIpN} Let $(X_i)_{i \in \mathbb Z}$ be a stationary
martingale differences sequence. Let $M_p=\bkE (|X_0|^p)$. Then for
any $p$ in $]2, 3]$ and any natural integer $N$,
\begin{equation} \label{boundIpN}
2^{-2N/p} \zeta_{p,N}^{2/p} \leq \Bigl( M_p +
\frac{1}{2 \sqrt 2} \sum_{K=0}^{N}2^{K(p/2-2)}\Vert Z_K
\Vert_{1, \Phi, p} \Bigr)^{2/p} + \frac{2}{p} \Delta_N
 \, ,
\end{equation}
where $Z_K = \bkE (S^2_{2^K}|{\mathcal F_0})-\bkE (S^2_{2^K})$ and
$\Delta_N=\sum_{K=0}^{N-1}2^{-2K/p}\Vert Z_K \Vert_{p/2}$.
\end{Proposition}
\noindent {\bf Proof of Proposition \ref{propIpN}.} The proof is
done by induction on $N$.  Let $\Y$ be a
sequence of $N(0, 1)$-distributed independent random
variables, independent of the sequence $\XZ$. For $m>0$, let  $T_m = Y_1 + Y_2 + \cdots + Y_{m}$.
Set $S_0 =T_0= 0$.  For $f$ numerical function  and $m\leq n$, set
$$
f_{n-m} (x) = \bkE ( f (x+ T_n-T_m) ) .
$$
Then, from the independence of the above sequences, \beq
\label{telessum} \bkE ( f ( S_n ) - f ( T_n ) ) = \sum_{ m=1}^n
D_{m} \ \hbox{ with }\ D_{m} = \bkE  \big( f_{n-m} (S_{m-1}   +  X_m
) - f_{n-m}  ( S_{m-1}  +  Y_m  ) \big) . \eeq Next, from the Taylor
integral formula at order two, for any  two-times differentiable
function $g$ and any $q$ in $]2, 3]$,
\begin{eqnarray*}
|g(x+h) - g(x) - g' (x) h - { \textstyle {1\over 2}} h^2 g'' (x)|
& \leq &
h^2 \int_0^1 (1-t) |g'' (x+th) - g'' (x)| dt \\
& \leq & h^2 \int_0^1 (1-t) |th|^{q-2} |g|_{\Lambda_q} dt ,
\end{eqnarray*}
whence \beq \label{taylorexp} | g(x+h) - g(x) - g' (x) h - {1\over
2} h^2 g'' (x)  | \leq {1 \over q(q-1)}  |h|^q |g|_{\Lambda_q} .
\eeq Let
$$
D'_m = \bkE ( f''_{n-m} (S_{m-1}) (X_m^2-1) ) = \bkE ( f''_{n-m}
(S_{m-1}) (X_m^2-Y_m^2) )
$$
From (\ref{taylorexp}) applied twice with $g= f_{n-m}$, $x =
S_{m-1}$ and $h=X_m$ or $h=Y_m$ together with the martingale property,
$$
\Big |D_m - \frac{1}{2} D'_m \Big |\leq   \frac{1}{p(p-1)}\,
|f_{n-m}|_{\Lambda_p} \bkE ( |X_{m}|^p + |Y_{m}|^p ) .
$$
Now $\bkE ( |Y_m|^p) \leq p-1 \leq (p-1) M_p$. Hence \beq
\label{reste}| D_m - (D'_m/2)|\leq M_p |f_{n-m}|_{\Lambda_p} \eeq
Moreover, if $f$ belongs to $\Lambda_p$, then the smoothed function
$f_{n-m}$ belongs to $\Lambda_p$. Hence, summing on $m$, we get that
\begin{equation}\label{summingn}
\bkE ( f (S_n ) - f ( T_n ) )  \leq  nM_p + ( D'/2)\quad \text{where
$D'=D'_1+ D'_2+\cdots + D'_n$}.
\end{equation}
\par\ssk
Suppose now that $n= 2^N$. To bound up  $D'$, we introduce a dyadic
scheme.
\begin{Notation}
Set $m_0 = m-1$ and write $m_0$ in basis $2$: $ m_0 = \sum_{i=0}^N
b_i 2^i $ with  $b_i=0$ or $b_i=1$ (note that $b_N  = 0$). Set $ m_L
=   \sum_{i=L}^N b_i 2^i , $ so that $m_N = 0$. Let $I_{L,k} =  ]k
2^L , (k+1) 2^L] \cap {\mathbb N}$ (note that $I_{N,1} = ]2^N,
2^{N+1}]$), $U_L^{(k)} = \sum_{i \in I_{L,k}  }  X_i$ and $\tilde
U_L^{(k)} = \sum_{i \in I_{L,k}  }  Y_i$. For the sake of brevity,
let $U_{L}^{(0)}=U_L$ and $\ti U_{L}^{(0)}= \ti U_L$.
\end{Notation}
Since  $m_N=0$, the following elementary identity is valid
$$
D'_m =  \sum_{L=0}^{N-1} \bkE \Bigl( ( f''_{n-1-m_L}  (S_{ m_L } ) -
f''_{n-1-m_{L+1} } (S_{m_{L+1}}  )) (X_m^2 -1) \Bigr) .
$$
Now $m_L \neq m_{L+1}$ only if $b_L = 1$, then in this case $m_L=k
2^L$ with $k$ odd. It follows that \beq D' = \sum_{L=0}^{N-1}\sum_{k
\in I_{N-L ,0} \atop k \text { odd }} {\mathbb E}\Big ((
f''_{n-1-k2^L} (S_{ k 2^L }  ) - f''_{n-1 - (k-1)2^L} (S_{(k-1) 2^L}
) )\sum_{\{m : m_L = k 2^L\}} (X_m^2 - \sigma^2)\Big) \, . \eeq Note
that $\{m : m_L = k 2^L\}=I_{L,k}$. Now by the martingale property
$$ {\mathbb E}_{k2^L}\Big(\sum_{ i\in I_{L,k}}
(X_i^2-\sigma^2)\Big)=  {\mathbb E}_{k 2^L} ( ( U_L^{(k)} )^2 ) -
\bkE ( ( U_L^{(k)} )^2 ):=Z_L^{(k)}\, .
$$
Since  $\X$ and $\Y$ are independent, we infer that
\begin{equation}\label{decdeltadiadique}
D' = \sum_{L=0}^{N-1}\sum_{k \in I_{N-L ,0} \atop k \text { odd }}
\bkE \Big ( \big ( f''_{n-1- k2^L}  (S_{ k 2^L } ) - f''_{n-1-k2^L}
(S_{(k-1) 2^L} + T_{k2^L} - T_{(k-1) 2^L}) \big ) Z_{L}^{(k)} \Big )
\, .
\end{equation}
By using (\ref{defgamrrio}), we get that
$$
D'  \leq  \sum_{L=0}^{N-1}\sum_{k \in I_{N-L ,0} \atop k \text { odd
}}  \bkE ( |U_{L}^{(k-1)} - \ti U_{L}^{(k-1)}|^{p-2} |Z_{L}^{(k)}| )
\, .
$$
From the stationarity of $\X$ and the above inequality,
\beq \label{bounddelta} D' \leq \frac{1}{2} \sum_{K=0}^{N-1} 2^{N-K}
\bkE ( |U_K - \ti U_K|^{p-2} |Z_{K}^{(1)}| ) . \eeq Now let $V_K$ be
the $N(0, 2^K )$-distributed random variable defined from $U_K$ via
the quantile transformation, that is
$$
V_K = 2^{K/2} \Phi^{-1} ( F_K ( U_K - 0 ) + \delta_K (F_K (
U_K ) - F_K (U_K - 0) ) )
$$
where  $F_K$ denotes  the d.f. of $U_K$, and $(\delta_K)$ is a
sequence of independent uniformly distributed r.v.'s, independent of
the underlying random variables. Now, from the subadditivity of $x \rightarrow x^{p-2}$,
$|U_K- \ti U_K|^{p-2} \leq |U_K - V_K|^{p-2} + |V_K - \ti U_K|^{p-2}$.
Hence
\beq
\label{majUZ1} \bkE ( |U_K - \ti U_K|^{p-2} |Z_{K}^{(1)}|
) \leq \Vert U_K - V_K \Vert_p^{p-2} \Vert Z_{K}^{(1)}\Vert_{p/2} +
\bkE ( |V_K - \ti U_K|^{p-2} |Z_{K}^{(1)}| ) \, .
\eeq
By definition
of $V_K$, the real $\Vert U_K - V_K \Vert_p$ is the so-called
Wasserstein distance of order $p$ between the law of $U_K^{(0)}$ and
the $N (0, 2^K )$ normal law. Therefrom, by Theorem 3.1 of Rio
(2007) (which improves the constants given in Theorem 1 of Rio
(1998)), we get that
\begin{equation} \label{riolemma}
\Vert U_K - V_K \Vert_p \leq 2 (2(p-1) \zeta_{p,K} )^{1/p} .
\end{equation}
Now, since  $V_K$ and  $\ti U_K$ are independent, their difference
has the $N (0, 2^{K+1})$ distribution. Hence, by definition
of the envelope norm $\Vert \, . \, \Vert_{1, \Phi, p}$,
\begin{equation} \label{boundquant}
\bkE   ( |V_K - \ti U_K|^{p-2} |Z_K^{(1)}| ) \leq
2^{(K+1)(p/2-1)} \Vert Z_K \Vert_{1, \Phi, p} .
\end{equation}
From (\ref{majUZ1}), (\ref{riolemma}) and (\ref{boundquant}), we get
that
\beq \label{majUZ2} \bkE ( |U_K - \ti U_K|^{p-2} |Z_{K}^{(1)}|
) \leq 2^{p-4/p} \zeta_{p,K}^{p-2\over p} \Vert Z_K \Vert_{p/2}+
2^{(K+1)(p/2-1)} \Vert Z_K \Vert_{1, \Phi, p} \, . \eeq
 Then, from
(\ref{summingn}), (\ref{bounddelta}) and (\ref{majUZ2}), we get
$$
2^{-N} \zeta_{p,N} \leq  M_p  + 2^{p/2-3}
 \Delta'_N + 2^{p-2-4/p}  \sum_{K=0}^{N-1} 2^{-K}
\zeta_{p,K}^{p-2\over p} \Vert Z_K \Vert_{p/2}    ,
$$
where $ \Delta'_N=\sum_{K=0}^{N-1}2^{K(p/2-2)}\Vert Z_K \Vert_{1,
\Phi, p}$. Consequently we get the induction
inequality
 \begin{equation}
\label{indineg} 2^{-N} \zeta_{p,N} \leq
M_p   + \frac{1}{2 \sqrt 2}\Delta'_N +
 \sum_{K=0}^{N-1} 2^{-K} \zeta_{p,K}^{p-2\over p} \Vert Z_K
\Vert_{p/2}  \, .
\end{equation}
We now prove (\ref{boundIpN}) by induction on $N$. Assume that
$\zeta_{p,L}$ satisfies (\ref{boundIpN}) for any $L$ in $[0,N-1]$.
Starting from (\ref{indineg}), using the induction hypothesis and
the fact that $\Delta'_K \leq \Delta'_N$, we  get that
$$
 2^{-N} \zeta_{p,N} \leq M_p    + \frac{1}{2\sqrt 2}\Delta'_N +
 \sum_{K=0}^{N-1} 2^{-2K/p}  \Vert Z_K \Vert_{p/2}
\Big ( \Big(  M_p + \frac{1}{2\sqrt 2}
\Delta'_N\Big )^{2/p} + \frac{2}{p} \Delta_K \Big)^{p/2-1} \, .
$$
Now $2^{-2K/p}\Vert Z_K \Vert_{p/2} = \Delta_{K+1} - \Delta_{K}$.
Consequently
\begin{eqnarray*}
2^{-N} \zeta_{p,N}  \leq M_p    + \frac{1}{2 \sqrt 2}\Delta'_N +
 \int_0^{\Delta_N} \Big ( \Big (M_p +
\frac{1}{2\sqrt 2}\Delta'_N \Big)^{2/p}+ \frac{2}{p} x \Big )^{p/2-1} dx \, ,
\end{eqnarray*}
which implies  (\ref{boundIpN})  for $\zeta_{p,N}$. $\square$
\par\medskip
In order to prove Theorem \ref{Thmart}, we will also need a smoothing argument. This is the purpose
of the lemma below.

\par\msk\no
\begin{Lemma} \label{smoothing} For any $r$ in $)0, p]$,
$\zeta_r ( P_{S_n} , G_n ) \leq 2 \zeta_r ( P_{S_n} * G_1 , G_n * G_1) + 4\sqrt{2}$.
\end{Lemma}
\noindent {\bf Proof of Lemma \ref{smoothing}.} Throughout the sequel, let $Y$
be a $N(0,1)$-distributed random variable, independent of the $\sigma$-field generated by the random variables
$(X_i)_i$ and $(Y_i)_i$.
\par\ssk
For $r\leq 2$, since $\zeta_r$ is an ideal metric with respect to the convolution,
$$
\zeta_r ( P_{S_n} , G_n ) \leq  \zeta_r ( P_{S_n} * G_1 , G_n * G_1) + 2 \zeta_r (\delta_0 , G_1) \leq
\zeta_r ( P_{S_n} * G_1 , G_n * G_1) + 2 \bkE |Y|^r
$$
which implies Lemma \ref{smoothing} for $r\leq 2$.
For $r>2$, from (\ref{taylorexp}), for any $f$ in $\Lambda_r$,
$$
\textstyle{
f (S_n) -f(S_n+Y) + f'(S_n) Y - {1\over 2} f'' (S_n) Y^2 \leq {1 \over r(r-1)} |Y|^r  . }
$$
Taking the expectation and noting that $\bkE |Y|^r \leq r-1$ for $r$ in $]2,3]$, we infer that
$$
\textstyle{
\bkE ( f (S_n) -f(S_n+Y) - {1\over 2} f'' (S_n) ) \leq {1 \over r} . }
$$
Obviously this inequality still holds for $T_n$ instead of $S_n$ and $-f$ instead of $f$, so that
adding the so obtained inequality,
$$
\textstyle{
\bkE ( f(S_n) - f (T_n)  \leq   \bkE ( f(S_n+Y) - f(T_n+Y) ) + {1\over 2} \bkE ( f'' (S_n) - f''(T_n)) + 1  . }
$$
It follows that
$$
\textstyle{
\zeta_r ( P_{S_n} , G_n ) \leq  \zeta_r ( P_{S_n} * G_1 , G_n * G_1) + {1\over 2} \zeta_{r-2} (P_{S_n} , G_n) + 1. }
$$
Now $r-2 \leq 1$. Hence
$$
\zeta_{r-2} (P_{S_n} , G_n) = W_{r-2} ( (P_{S_n} , G_n) \leq ( W_r (P_{S_n} , G_n) )^{r-2} .
$$
Next, by Theorem 3.1 in Rio (2007),
$W_r (P_{S_n} , G_n) \leq ( 32 \zeta_r  (P_{S_n} , G_n) )^{1/r}$.  Furthermore
$$
(32 \zeta_r (P_{S_n} , G_n) )^{1-2/r} \leq \zeta_r (P_{S_n} , G_n)
$$
as soon as $\zeta_r (P_{S_n} , G_n) \geq 2^{(5r/2)-5}$. This condition holds for any $r$ in $]2,3]$ if
$\zeta_r (P_{S_n} , G_n) \geq 4\sqrt{2}$. Then, from the above inequalities
$$
\textstyle{
\zeta_r ( P_{S_n} , G_n ) \leq  \zeta_r ( P_{S_n} * G_1 , G_n * G_1) + {1\over 2} \zeta_r (P_{S_n} , G_n) + 1, }
$$
which implies Lemma \ref{smoothing}  $\square$
\par\msk

\par\msk
We go back to the proof of Theorem \ref{Thmart}. We will first
complete the proof in the case $p=r$. Next we will derive the
general case at the end of the proof.

Let $\zeta^*_{p,N} = \sup_{n \leq 2^N} \zeta_p( P_{S_n}, G_n)$. We
will bound up $\zeta^*_{p,N}$ by induction on $N$. Let $n \in ]2^N,
2^{N+1}]$. Hence $n = 2^N + \ell$ with $\ell \in [1, 2^N]$. We first
notice that \begin{equation*}  \zeta_r ( P_{S_n} , G_n ) \leq
\zeta_r ( P_{S_n} , P_{S_{\ell}}* G_{2^N} ) + \zeta_r (
P_{S_{\ell}}* G_{2^N} , G_{\ell}* G_{2^N} ) \, . \end{equation*}
Now, with the same notation as in the proof of Proposition
\ref{propIpN}, we have
$$
\zeta_r ( P_{S_{\ell}}* G_{2^N} , G_{\ell}* G_{2^N} ) = \sup_{f \in
{\Lambda_r}} {\mathbb E} ( f_{2^N} (S_{\ell}) - f_{2^N} (T_{\ell}))
\leq |f * \phi_{2^{N/2}} |_{\Lambda_p} \zeta_p ( P_{S_{\ell}} ,
G_{\ell} ) \, .
$$
Applying Lemma \ref{regul}, we infer that \beq \label{bII}  \zeta_r
( P_{S_n} , G_n ) \leq \zeta_r ( P_{S_n} , P_{S_{\ell}}* G_{2^N}
)+c_{r,p} 2^{N(r-p)/2} \zeta_p ( P_{S_{\ell}} , G_{\ell} ) \, . \eeq
On the other hand, setting $\tilde S_{\ell} = X_{1-\ell} + \cdots +
X_0$, we have that $S_n$ is distributed as $ \tilde S_{\ell} +
S_{2^N}$. Using Lemma \ref{smoothing}, we then derive that
 \beq \label{bI}
\zeta_r ( P_{S_n} , P_{S_{\ell}}* G_{2^N} ) \leq 4 \sqrt{2} +
2\sup_{f \in {\Lambda_r}} {\mathbb E} ( f (\tilde S_{\ell} + S_{2^N}
+ Y) - f (\tilde S_{\ell} + T_{2^N}+ Y))
 \eeq
Let $D'_m= {\mathbb E}  ( f''_{2^N -m+1} (\tilde S_{\ell} + S_{m-1})
(X_m^2 -1) )$. Following the proof of Proposition \ref{propIpN}, we
get that
\begin{equation} \label{b2I}
 {\mathbb E} ( f (\tilde S_{\ell} + S_{2^N} + Y) - f (\tilde
S_{\ell} + T_{2^N}+ Y))  = ( D'_1 + \cdots + D'_{2^N} )/2 + R_1 + \cdots + R_{2^N}  \, ,
\end{equation}
where, as in (\ref{reste}),
\begin{equation}\label{reste2}
R_m \leq M_p |f_{2^N-m+1}|_{\Lambda_p}\, .
\end{equation}
In the case $r=p-2$, we will need the more precise upper bound
\begin{equation}
R_m  \leq  \bkE \Big ( X_m^2 \big( \Vert f''_{2^N -m+1}
\Vert_\infty \wedge { \frac{1}{6} }\Vert f_{2^N -m+1} ^{(3)}
\Vert_\infty |X_m| \big ) \Big )  +
 {\frac{1}{6} }\Vert f_{2^N -m+1} ^{(3)} \Vert_\infty \bkE (|Y_m|^3)
 \, ,
\end{equation}
which is derived from the Taylor formula at orders two and three.
From (\ref{reste2}) and Lemma \ref{regul}, we have that
\begin{equation} \label{borner1}
R : =  R_1 + \cdots + R_{2^N} = O(2^{N(r-p+2)/2})\quad \text{if $r> p -2$,
and}\quad R =O(N) \quad \text{if $(r,p)=(1,3)$}\,
.
\end{equation}

It remains to consider the case $r=p-2$ and $r<1$. Applying Lemma
\ref{regul}, we get that for $i \geq 2$,
\begin{equation}\label{bornefi}
\Vert
f^{(i)}_{2^N -m+1} \Vert_\infty \leq c_{r,i} (2^N -m +1)^{(r-i)/2}
\, .
\end{equation}
It follows that
\begin{eqnarray*}
 \sum_{m=1}^{2^N}\bkE \Big ( X_{m}^2 \big( \Vert f''_{2^N -m+1}
\Vert_\infty \wedge \Vert f_{2^N -m+1} ^{(3)} \Vert_\infty |X_{m}|
\big ) \Big ) &\leq& C \sum_{m=1}^{\infty}
\frac{1}{m^{1-r/2}}{\mathbb E} \left (
X_0^2 \Bigl( 1 \wedge \frac{|X_0|}{\sqrt m} \Bigr) \right ) \\
 &\leq & C{\mathbb E} \Big ( \sum_{m=1}^{[X_0^2]}\frac{X_0^2
}{m^{1-r/2}} +\sum_{m=[X_0^2] +
1}^{\infty}\frac{|X_0|^3}{m^{(3-r)/2}} \Big )\, .
\end{eqnarray*}
Consequently for $r=p-2$ and $r<1$,
\begin{equation} \label{borner2}
R_1 + \cdots + R_{2^N}  \leq C(M_p+{\mathbb E}(|Y|^3))\, .
\end{equation}

We  now bound up $D'_1 + \cdots + D'_{2^N} $. Using the dyadic scheme as
in the proof of Proposition \ref{propIpN}, we get that
\begin{eqnarray*} D'_m &=& \sum_{L=0}^{N-1} \bkE \Bigl( (
f''_{2^N-m_L}  (\tilde S_{\ell}  + S_{ m_L } ) - f''_{2^N-m_{L+1} }
(\tilde S_{\ell} + S_{m_{L+1}}  ) (X_m^2 -1) \Bigr) +  \bkE  (
f''_{2^N} (\tilde S_{\ell}  ) (X_m^2 -1
) ) \\
& := &  D''_m + \bkE  ( f''_{2^N} (\tilde S_{\ell}  ) (X_m^2 -1 ) )
\, .
\end{eqnarray*}

Notice first that
\begin{equation*}
\sum_{m=1}^{2^N}\bkE  ( f''_{2^N} (\tilde S_{\ell}  ) (X_m^2 -1 ) )
= \bkE ( (f''_{2^N} (\tilde S_{\ell}  )- f''_{2^N} (T_{\ell}
))Z_N^{(0)}) \, .
\end{equation*}
Hence using Lemma \ref{regul}, we get that
\begin{equation*}
\sum_{m=1}^{2^N}\bkE  ( f''_{2^N} (\tilde S_{\ell}  ) (X_m^2 -1 ) )
\leq C 2^{N(r-p)/2} \bkE  ( |\tilde S_{\ell} -  T_{\ell}
|^{p-2}|Z_N^{(0)} | )
  \, .
\end{equation*}
Proceeding as to get (\ref{majUZ2}), we have that
\begin{equation*} \bkE  ( |\tilde S_{\ell} -  T_{\ell}
|^{p-2}|Z_N^{(0)} |)  \leq  2^{p-4/p}(\zeta_p(P_{S_{\ell}} ,
G_{\ell}) )^{(p-2)/p}\| Z_N^{(0)} \|_{p/2} + (2 \ell)^{p/2 -1} \|
Z_N^{(0)} \|_{1, \Phi, p} \, .
\end{equation*}
Using Remark \ref{rmkequivalence},  (\ref{C1}) and (\ref{C2}) entail
that $\| Z_N^{(0)} \|_{p/2} = o( 2^{2N/p})$ and $\| Z_N^{(0)} \|_{1,
\Phi, p} =o(2^{N(2-p/2})$. Hence, for some $\epsilon(N)$ tending to
$0$ as $N$ tends to infinity, one has
\begin{equation} \label{boundtermsup}
 D'_1 + \cdots + D'_{2^N}  \leq  C(\epsilon(N) 2^{N((r-p)/2+2/p} (\zeta_p(P_{S_{\ell}} ,
G_{\ell}) )^{(p-2)/p} +  2^{N(r+2-p)/2}) \, .
\end{equation}

Next, proceeding as in the proof of (\ref{decdeltadiadique}), we get
that
$$
\sum_{m=1}^{2^N}D''_m \leq \sum_{L=0}^{N-1}\sum_{k \in I_{N-L ,0}
\atop k \text { odd }} \bkE \Big ( \big ( f''_{2^N- k2^L} (\tilde
S_{\ell} + S_{ k 2^L } ) - f''_{2^N-k2^L} (\tilde S_{\ell} +
S_{(k-1) 2^L} + T_{k2^L} - T_{(k-1) 2^L}) \big ) Z_{L}^{(k)} \Big )
\, .
$$
If $r>p-2$ or $(r,p)=(1,3)$, from Lemma \ref{regul}, the
stationarity of $\X$ and the above inequality,
\begin{eqnarray*}\sum_{m=1}^{2^N}D''_m  & \leq & C
\sum_{L=0}^{N-1}\sum_{k \in I_{N-L ,0} \atop k \text { odd }}
(2^N-k2^L)^{(r-p)/2}\bkE  \big (  | U_{L} - \tilde U_{L}|^{p-2} \big
| Z_{L}^{(1)} \big | \big) \, .
\end{eqnarray*}
It follows that
\begin{eqnarray}\label{b1D''}
\sum_{m=1}^{2^N}D''_m &\leq& C {2^N}^{(r+2-p)/2}\sum_{L=0}^{N}
2^{-L}\mathbb E \big ( \big | U_{L}- \tilde U_{L} \big|^{p-2} \big |
Z_{L}^{(1)} \big | \big ) \quad \text{if $r>p-2$,}\\
\label{b2D''} \sum_{m=1}^{2^N}D''_m&\leq &C N \sum_{L=0}^{N}
2^{-L}\mathbb E \big ( \big | U_{L}- \tilde U_{L} \big| \big |
Z_{L}^{(1)} \big | \big ) \quad \text{if $r=1$ and $p=3$.}
\end{eqnarray}

In the case $r=p-2$ and $r>1$, we have
\begin{eqnarray*}\sum_{m=1}^{2^N}D''_m  & \leq & C
\sum_{L=0}^{N-1}\sum_{k \in I_{N-L ,0} \atop k \text { odd }} \bkE
\Big ( \big ( \|f''_{2^N -k2^L}\|_{\infty} \wedge \|f'''_{2^N -k2^L
}\|_{\infty} \big | U_{L} - \tilde U_{L} \big| \big ) \big |
Z_{L}^{(1)} \big | \Big ) \, .
\end{eqnarray*}
Applying (\ref{bornefi}) to $i=2$ and $i=3$, we obtain
$$ \sum_{m=1}^{2^N}D''_m   \leq  C \sum_{L=0}^{N} 2^{(r-2)L/2}\mathbb E \Big (\big | Z_{L}^{(1)}\big |\sum_{k=1}^{2^{N-L}}
k^{(r-2)/2}\big ( 1 \wedge \frac{1}{2^{L/2}\sqrt k}\big | U_{L}-
\tilde U_{L} \big| \big )  \Big ) \, ,
 $$
 Proceeding as to get (\ref{borner2}), we have that
 $$
\sum_{k=1}^{2^{N-L}} k^{(r-2)/2}\big ( 1 \wedge
\frac{1}{2^{L/2}\sqrt k}\big | U_{L}- \tilde U_{L} \big| \big ) \leq
\sum_{k=1}^{\infty} k^{(r-2)/2}\big ( 1 \wedge \frac{1}{2^{L/2}\sqrt
k}\big | U_{L}- \tilde U_{L} \big| \big ) \leq C| U_{L}- \tilde
U_{L} \big|^r \, .
 $$
 It follows that
 \begin{equation}\label{b3D''}
\sum_{m=1}^{2^N}D''_m   \leq  C \sum_{L=0}^{N} 2^{-L}\mathbb E \Big
( \big | U_{L}- \tilde U_{L} \big|^{r} \big |Z_{L}^{(1)} \big | \Big
)\quad \text{if  $r = p-2$ and $r<1$.}
 \end{equation}

 Now by Remark \ref{rmkequivalence}, (\ref{C1}) and (\ref{C2}) are
respectively equivalent to
$$
\sum_{K\geq 0}2^{K(p/2-2)}\Vert Z_K \Vert_{1, \Phi, p} < \infty \,
,\, \text{ and } \, \sum_{K \geq 0}2^{-2K/p}\Vert Z_K \Vert_{p/2} <
\infty \, .$$
 Next, by Proposition \ref{propIpN}, $\zeta_{p,K}=O(2^K)$
under (\ref{C1}) and (\ref{C2}). Therefrom, taking into account the
inequality (\ref{majUZ2}), we derive that under (\ref{C1}) and
(\ref{C2}), \begin{eqnarray} \label{importantbound} 2^{-L}\bkE \Big
( \big | U_{L} - \tilde U_{L} \big|^{p-2} \big | Z_{L}^{(1)} \big |
\Big  ) & \leq & C 2^{-2L/p}\Vert Z_L \Vert_{p/2} + C 2^{L(p/2-2)}
\Vert Z_K \Vert_{1, \Phi, p} \, .
\end{eqnarray}
Consequently, combining (\ref{importantbound}) with the upper bounds
(\ref{b1D''}), (\ref{b2D''}) and (\ref{b3D''}), we obtain that
\begin{equation}\label{dernierD''}
\sum_{m=1}^{2^N}D''_m= \left\{
\begin{array}{ll} O(2^{N(r+2-p)/2})& \text{ if $r \geq p-2$ and $(r,p) \neq (1,3)$} \\
O( N ) & \text{ if $r=1$ and $p=3$.}
\end{array}
\right. \end{equation} From (\ref{bII}), (\ref{bI}), (\ref{b2I}),
(\ref{borner1}), (\ref{borner2}), (\ref{boundtermsup}) and
(\ref{dernierD''}), we get that if $r \geq p-2$ and $(r,p) \neq
(1,3)$, \beq \label{bound1*} \zeta_r ( P_{S_n} , G_n ) \leq
c_{r,p}2^{N(r-p)/2} \zeta_p ( P_{S_{\ell}} , G_{\ell} )+
C(2^{N(r+2-p)/2}+ 2^{N( (r-p)/2 +2/p)}\epsilon(N)
(\zeta_p(P_{S_{\ell}} , G_{\ell}) )^{(p-2)/p}  ) \eeq and if $r=1$
and $p=3$, \beq \label{bound2*} \zeta_1 ( P_{S_n} , G_n ) \leq  C (N
+ 2^{-N} \zeta_3 ( P_{S_{\ell}} , G_{\ell} )+ 2^{-N/3}
(\zeta_3(P_{S_{\ell}} , G_{\ell}) )^{1/3} ) \, . \eeq Since
$\zeta^*_{p,N} = \sup_{n \leq 2^N} \zeta_p( P_{S_n}, G_n)$, we infer
from (\ref{bound1*}) applied to $r=p$ that
$$
\zeta^*_{p,N+1}  \leq   \zeta^*_{p,N}  + C(2^{N}+
2^{2N/p}\epsilon(N) (\zeta^*_{p,N})^{(p-2)/p} ) \, .$$ Let $N_0$ be
 such that $C\epsilon(N) \leq 1/2$ for $N\leq N_0$, and let $K\geq 1$ be such that
 $\zeta^*_{p, N_0}\leq K 2^{N_0}$. Choosing $K$ large enough such
 that $K\geq 2C$, we can easily prove by induction that $\zeta^*_{p,
 N} \leq K 2^N$ for any $N\geq N_0$. Hence Theorem \ref{Thmart} is proved in the case
 $r=p$.

 For $r$ in $[p-2, p[$, Theorem \ref{Thmart} follows by taking into
 account the bound $\zeta^*_{p,N} \leq K 2^N$, valid  for any $N\geq N_0$, in the inequalities
 (\ref{bound1*}) and (\ref{bound2*}).

\subsection{Proof of Theorem \ref{propNAS}} By
(\ref{Cond1cob}), we get that (see Voln\'y (1993))
\begin{equation}
\label{cobound} X_0=D_0 + Z_0-Z_0 \circ T,
\end{equation} where $$Z_0=\sum_{k=0}^\infty {\mathbb E}(X_k|{\mathcal F}_{-1})-
\sum_{k=1}^\infty (X_{-k}-{\mathbb E}(X_{-k}|{\mathcal
F}_{-1})) \quad \text{and} \quad D_0=\sum_{k \in {\mathbb Z}}
{\mathbb E}(X_k|{\mathcal F}_0)-{\mathbb E}(X_k|{\mathcal F}_{-1})\,
.
$$
Note that $D_0 \in {\mathbb L}^p$, $D_0 $ is ${\mathcal
F}_0$-measurable, and $\bkE (D_0 |{\mathcal F}_{-1})=0$. Let
$D_i=D_0 \circ T^i$, and $Z_i=Z_0 \circ T^i$. We obtain that
\begin{equation}
\label{coboundSm} S_{n}=M_{n} + Z_1-Z_{n+1} \, ,
\end{equation}
where $M_{n}=\sum_{j=1}^nD_j$.  We first bound up ${\mathbb E} ( f
(S_{n}) - f(M_{n}))$ by using the following lemma
\begin{Lemma} \label{compSnMn} Let $p
\in ]2,3]$ and $r \in [p-2,p]$. Let $(X_i)_{i \in \mathbb Z}$ be a
stationary sequence of centered random variables in ${\mathbb
L}^{2\vee r}$. Assume that $S_n = M_n + R_n$ where
$(M_n-M_{n-1})_{n>1}$ is a strictly stationary sequence of
martingale differences in $\mathbb L^{2\vee r}$, and $R_n$ is such
that $\bkE(R_n)=0$. Let
$n\sigma^2= {\mathbb E}(M_n^2)$, $n\sigma_n^2={\mathbb E}(S_n^2)$
and $\alpha_n=\sigma_n/\sigma$.
\begin{enumerate}
\item If $r \in [p-2,1]$ and $\bkE|R_n|^r=O(n^{(r+2-p)/2})$, then
$\displaystyle \zeta_r(P_{S_n}, P_{M_n}) =O(n^{(r+2-p)/2}).$
 \item If $r \in ]1,2]$ and $\|R_n\|_r=O(n^{(3-p)/2})$, then
$\displaystyle \zeta_r(P_{S_n}, P_{M_n}) = O (n^{(r+2-p)/2}).$
\item
If $r \in ]2,p]$, $\sigma^2>0$ and $\|R_n\|_r=O(n^{(3-p)/2})$, then
$\displaystyle \zeta_r(P_{S_n}, P_{\alpha_nM_n}) = O
(n^{(r+2-p)/2}).$
\item If
$r \in ]2,p]$, $\sigma^2=0$ and $\|R_n\|_r=O(n^{(r+2-p)/2r})$, then
$\displaystyle \zeta_r(P_{S_n}, G_{n \sigma_n^2}) = O
(n^{(r+2-p)/2}).$
\end{enumerate}
\end{Lemma}
\begin{Remark}
All the assumptions of  Lemma \ref{compSnMn} are satisfied as soon
as $\sup_{n >0} \|R_n\|_p < \infty$.
\end{Remark}
\noindent {\bf Proof of Lemma \ref{compSnMn}}. For $r \in ]0,1]$,
$\zeta_r(P_{S_n}, P_{M_n})\leq \bkE (| R_n |^r)$, which implies item
1. If $f \in \Lambda_r$ with $r \in ]1,2]$, from the Taylor integral
formula and since $\bkE ( R_n) = 0$, we get
\begin{eqnarray*}
{\mathbb E} ( f (S_{n}) - f(M_{n})) & = & {\mathbb
E}\Big(R_n\Big(f'(M_{n}) - f'(0) +
\int_0^1 (f'(M_{n} +t(R_n))-f'(M_{n}))dt \Big ) \Big)\\
& \leq & \Vert R_n \Vert_r \Vert f'(M_{n})- f'(0) \Vert_{r/(r-1)} +
\Vert R_n \Vert_r^r   \leq  \Vert R_n \Vert_r \Vert M_{n}
\Vert_{r}^{r-1} + \Vert R_n \Vert_r^r\, .
\end{eqnarray*} Since  $\|M_n\|_r\leq \|M_n\|_2 = \sqrt n \sigma$, we infer that
$\zeta_r(P_{S_n}, P_{M_n})=O (n^{(r+2-p)/2})$.

 Now if $f \in \Lambda_r$ with $r \in ]2,p]$ and if $\sigma>0$, we
define $g$ by $$g(t) = f(t) - t f'(0) - t^2 f''(0)/2\, .$$ The
function $g$ is then also in $ \Lambda_r$ and is such that
$g'(0)=g''(0)=0$. Since $ \displaystyle \alpha^2_{n} {\mathbb E
(M_{n}^2)}= {\mathbb E (S_{n}^2)}$, we have \beq \label{eq} {\mathbb
E} ( f (S_{n}) - f(\alpha_{n} M_{n})) = {\mathbb E} ( g (S_{n}) - g(
\alpha_{n} M_{n})) \, . \eeq Now from the Taylor integral formula at
order two, setting $\tilde R_{n} = R_n + (1- \alpha_{n}) M_{n}$,
\begin{eqnarray*}
{\mathbb E} ( g (S_{n}) - g(\alpha_n M_{n})) & = & {\mathbb
E}(\tilde R_n g'(\alpha_{n} M_{n}) )+ \frac 12 {\mathbb
E}((\tilde R_{n})^2 g''(\alpha_{n} M_{n}) ) \\
& & \quad + {\mathbb E}\big ((\tilde R_{n})^2
\int_0^1 (1-t)(g''(\alpha_{n} M_{n} +t\tilde R_{n})-g''(\alpha_{n} M_{n}))dt \big ) \\
& \leq & \frac{1}{r-1}{\mathbb E}(|\tilde R_{n}||\alpha_n
M_{n}|^{r-1} ) + \frac 12 \Vert \tilde R_{n} \Vert^2_r \Vert
g''(\alpha_{n} M_{n})\Vert_{r/(r-2)} + \frac 12
\Vert \tilde R_{n} \Vert_r^r \\
& \leq & \frac{1}{r-1} \alpha_{n}^{r-1}\Vert \ti R_{n} \Vert_r
\Vert M_{n} \Vert_{r}^{r-1} + \frac 12 \alpha_{n}^{r-2} \Vert \ti
R_{n} \Vert^2_r \Vert M_{n} \Vert_{r}^{r-2} + \Vert \ti R_{n}
\Vert_r^r\, .
\end{eqnarray*}
Now $\alpha_{n} = O(1)$ and $ \Vert \ti R_{n} \Vert_r \leq \Vert R_n
\Vert_r + |1 - \alpha_{n}| \Vert M_{n}\Vert_r$. Since
$|\|S_n\|_2-\|M_n\|_2|\leq  \|R_n\|_2$, we infer that
 $|1 -
\alpha_{n}| = O(n^{(2-p)/2})$. Hence, applying Burkh\"older's
inequality for martingales, we infer that $ \|\tilde R_n\|_r=
O(n^{(3-p)/2})$, and consequently $\zeta_r(P_{S_n}, P_{\alpha_n
M_n})=O (n^{(r+2-p)/2})$.

If $\sigma^2=0$, then $S_n=R_n$. Using that ${\mathbb E} ( f (S_{n})
- f(\sqrt n \sigma_n Y)) = {\mathbb E} ( g (R_{n}) - g(\sqrt n
\sigma_n Y))$, and applying again Taylor's formula, we obtain that
$$
\sup_{f \in \Lambda_r}|{\mathbb E} ( f (S_{n}) - f( \sqrt n \sigma_n
Y))| \leq \frac{1}{r-1} \Vert \bar R_{n} \Vert_r \Vert \sqrt n
\sigma_n Y \Vert_{r}^{r-1} + \frac 12  \Vert \bar R_{n} \Vert^2_r
\Vert \sqrt n \sigma_n Y \Vert_{r}^{r-2} + \Vert \bar R_{n}
\Vert_r^r\, ,
$$
where $\bar R_n= R_n -\sqrt n \sigma_n Y$. Since $\sqrt n \sigma_n
=\|R_n\|_2=O(n^{(r+2-p)/2r})$, the result follows. $\square$

\medskip
By (\ref{coboundSm}), we can apply Lemma \ref{compSnMn} with $R_n
:=Z_1 - Z_{n+1}$. Then for $ p-2 \leq  r \leq 2$, the result follows
if we prove that under (\ref{Cond2cob}), $M_n$ satisfies the
conclusion of Theorem \ref{Thmart}. Now if $2 < r \leq p$ and
$\sigma^2>0$, we first notice that
$$
\zeta_r(P_{\alpha_n M_n},G_{n
\sigma_n^2})=\alpha_n^r\zeta_r(P_{M_n},G_{n \sigma^2} )\, .
$$
Since $\alpha_{n} = O(1)$, the result will follow by Item 3 of Lemma
\ref{compSnMn}, if we prove that under (\ref{Cond2cob}), $M_n$
satisfies the conclusion of Theorem \ref{Thmart}. We shall prove
that \beq \label{condMartcob} \sum_{n \geq
1}\frac{1}{n^{3-p/2}}\Vert \bkE (M_n^2 | {\mathcal F}_0) - \bkE
(M_n^2) \Vert_{p/2} < \infty \, .\eeq In this way, both (\ref{C1})
and (\ref{C2}) will be satisfied. Suppose that we can show that \beq
\label{apMnSnp} \sum_{n \geq 1}\frac{1}{n^{3-p/2}}\Vert \bkE (M_n^2
| {\mathcal F}_0) - \bkE (S_n^2 | {\mathcal F}_0)\Vert_{p/2} <
\infty \, ,\eeq then by taking into account the condition
(\ref{Cond2cob}), (\ref{condMartcob}) will follow. Indeed, it
suffices to notice that (\ref{apMnSnp}) also entails that \beq
\label{cons1apMnSnp} \sum_{n \geq 1}\frac{1}{n^{3-p/2}} | \bkE
(S_n^2) - \bkE (M_n^2) | < \infty \, ,\eeq and to write that
\begin{eqnarray*}
\Vert \bkE (M_n^2 | {\mathcal F}_0) - \bkE (M_n^2 )\Vert_{p/2} &
\leq  &  \Vert \bkE (M_n^2 | {\mathcal F}_0) - \bkE (S_n^2 |
{\mathcal F}_0)\Vert_{p/2} \\& & \quad \quad +  \Vert \bkE (S_n^2 |
{\mathcal F}_0) - \bkE (S_n^2 )\Vert_{p/2} + |  \bkE (S_n^2 ) - \bkE
(M_n^2 ) |\, .
\end{eqnarray*} Hence, it remains to prove (\ref{apMnSnp}).  Since $S_n = M_n +
Z_1-Z_{n+1}$, and since $Z_i=Z_0\circ T^i$ is in ${\mathbb L}^p$,
(\ref{apMnSnp}) will be satisfied  provided that \beq
\label{conddouble} \sum_{n \geq 1}\frac{1}{n^{3-p/2}}\Vert
S_n(Z_1-Z_{n+1}) \Vert_{p/2} < \infty \, . \eeq Notice that
\begin{eqnarray*}
\Vert S_n(Z_1-Z_{n+1}) \Vert_{p/2} & \leq & \Vert M_n \Vert_{p}
\Vert Z_1-Z_{n+1} \Vert_{p} + \Vert Z_1-Z_{n+1} \Vert^2_{p}\, .
\end{eqnarray*}
From Burkholder's inequality, $\Vert M_n \Vert_p =O(\sqrt n)$ and
from (\ref{Cond1cob}), $\sup_n\Vert Z_1-Z_{n+1} \Vert_{p} < \infty$.
Consequently (\ref{conddouble}) is satisfied for any $p$ in $]2,3[$.

\medskip

\subsection{Proof of Theorem \ref{propNASp3}}  Starting from
(\ref{coboundSm}) we have that \beq \label{decM_n} M_n  :=  S_n +R_n
+ \tilde R_n \, ,\eeq where
$$ R_n = \sum_{k \geq n +1 } \bkE (X_k |{\mathcal F_n}) -\sum_{k
\geq 1 } \bkE (X_k |{\mathcal F_0}) \text{ and }\tilde R_n =
\sum_{k \geq 0 } (X_{-k}-\bkE (X_{-k} |{\mathcal F_0})) - \sum_{k
\geq -n } (X_{-k}-\bkE (X_{-k} |{\mathcal F_n})) \, .
$$Arguing as in the proof of Proposition \ref{propNAS} the proposition will follow from (\ref{Cond2cobp3}), if we prove that
\beq \label{apMnSn3} \sum_{n \geq 1}^{\infty}\frac{1}{n^{3/2}}\Vert
\bkE (M_n^2 | {\mathcal F}_0) - \bkE (S_n^2 | {\mathcal F}_0)
\Vert_{3/2} < \infty \, . \eeq Under (\ref{Condcobp3adap}), $\sup_{n
\geq 1}\|R_n \|_3 < \infty$ and $\sup_{n \geq 1} \|\Tilde R_n \|_3 <
\infty$. Hence (\ref{apMnSn3}) will be verified as soon as
 \beq \label{cond1doubleproduit} \sum_{n=1}^{\infty}\frac{1}{n^{3/2}}\Vert \bkE (
S_n (R_n + \Tilde R_n ) |{\mathcal F}_0)  \Vert_{3/2} < \infty \, .
\eeq We first notice that the decomposition (\ref{decM_n}) together
with Burkholder's inequality for martingales and the fact that
$\sup_n\|R_n \|_3 < \infty$ and $\sup_n \| \Tilde R_n  \|_3 <
\infty$, implies that \beq \label{majnormeS_n} \Vert S_n \Vert_3
\leq C \sqrt n \, .\eeq Now to prove (\ref{cond1doubleproduit}), we
first notice that
\begin{eqnarray} \label{pradap1}
& & \Big \Vert \bkE \Big ( S_n\sum_{k \geq 1 } \bkE (X_k |{\mathcal
F_0}) \Big| {{\mathcal F_0}} \Big )  \Big\Vert_{3/2} \leq
 \Vert \bkE  ( S_n | {{\mathcal F_0}}  )\Vert_3 \Big \Vert \sum_{k \geq 1 } \bkE (X_k |{\mathcal
F_0})\Big \Vert_{3} \, ,
\end{eqnarray}
which is bounded by using (\ref{Condcobp3adap}). Now write
\begin{eqnarray*}
\bkE \Big ( S_n \sum_{k \geq n +1 } \bkE (X_k |{\mathcal F_n}) \Big
| {{\mathcal F_0}} \Big )  =  \bkE \Big ( S_n \sum_{k \geq 2n+1}
\bkE (X_k |{\mathcal F_n}) \Big | {{\mathcal F_0}} \Big ) + \bkE  (
S_n \bkE (S_{2n}-S_n |{\mathcal F_n})  | {{\mathcal F_0}} ) \, .
\end{eqnarray*}
Clearly
\begin{eqnarray} \label{pradap2}
\Big \Vert \bkE \Big ( S_n \sum_{k \geq 2n+1} \bkE (X_k |{\mathcal
F_n}) \Big | {{\mathcal F_0}} \Big ) \Big\Vert_{3/2} & \leq & \Vert
S_n \Vert_{3} \Big \Vert \sum_{k \geq 2n+1 } \bkE
(X_k |{\mathcal F_n}) \Big \Vert_3 \nonumber \\
& \leq & C \sqrt n \Big \Vert \sum_{k \geq n+1 } \bkE (X_k
|{\mathcal F_0}) \Big \Vert_3 \, ,
\end{eqnarray}
by using (\ref{majnormeS_n}). Considering the bounds
(\ref{pradap1}) and (\ref{pradap2}) and the condition
(\ref{Condcobp3adap}), in order to prove that
 \beq \label{cond1doubleproduitadap} \sum_{n=1}^{\infty}\frac{1}{n^{3/2}}\Vert \bkE (
S_n R_n  |{\mathcal F}_0) \Vert_{3/2} < \infty \, ,\eeq it is
sufficient to prove that
 \beq \label{resteadap} \sum_{n=1}^{\infty}\frac{1}{n^{3/2}}\Vert \bkE  ( S_n \bkE (S_{2n}-S_n
|{\mathcal F_n}) | {{\mathcal F_0}}  ) \Vert_{3/2}< \infty \, .\eeq
With this aim,  take $p_n =[\sqrt n]$ and write
\begin{eqnarray} \label{decSnS2npn}
\bkE ( S_n \bkE (S_{2n}-S_n |{\mathcal F_n}) | {{\mathcal F_0}}  ) &
= & \bkE  ( (S_n - S_{n - p_n}) \bkE (S_{2n}-S_n |{\mathcal F_n}) |
{{\mathcal F_0}} ) \nonumber
\\& & \quad \quad + \bkE  ( S_{n - p_n} \bkE (S_{2n}-S_n |{\mathcal F_n}) |
{{\mathcal F_0}}  ) .
\end{eqnarray}
By stationarity and (\ref{majnormeS_n}), we get that
$$
\sum_{n=1}^{\infty}\frac{1}{n^{3/2}}\Vert \bkE  ( (S_n - S_{n -
p_n}) \bkE (S_{2n}-S_n |{\mathcal F_n}) | {{\mathcal F_0}} )
\Vert_{3/2} \leq C
\sum_{n=1}^{\infty}\frac{\sqrt{p_n}}{n^{3/2}}\Vert \bkE (S_n
|{\mathcal F_0}) \Vert_{3} \, ,
$$
which is finite under (\ref{Condcobp3adap}), since $p_n =[\sqrt n]$.
Hence from (\ref{decSnS2npn}), (\ref{resteadap}) will follow if we
prove that \beq \label{resteadapSpn}
\sum_{n=1}^{\infty}\frac{1}{n^{3/2}}\Vert \bkE  ( S_{n - p_n} \bkE
(S_{2n}-S_n |{\mathcal F_n}) | {{\mathcal F_0}}  )\Vert_{3/2}<
\infty \, .\eeq With this aim we first notice that
\begin{multline*}
 \Vert \bkE  ( (S_{n - p_n} - \bkE (S_{n - p_n}|{\mathcal
F_{n-p_n}}) \bkE (S_{2n}-S_n |{\mathcal F_n}) | {{\mathcal F_0}}
 )\Vert_{3/2} \\
 \leq \Vert S_{n - p_n} - \bkE (S_{n - p_n}|{\mathcal F_{n-p_n}})
\Vert_{3} \Vert \bkE (S_{2n}-S_n |{\mathcal F_n}) \Vert_{3} \,
,\end{multline*} which is bounded under (\ref{Condcobp3adap}).
Consequently (\ref{resteadapSpn}) will hold if we prove that \beq
\label{resteadapSpn2} \sum_{n=1}^{\infty}\frac{1}{n^{3/2}}\Vert \bkE
 ( \bkE (S_{n - p_n}|{\mathcal F_{n-p_n}})\bkE (S_{2n}-S_n
|{\mathcal F_n}) | {{\mathcal F_0}} )\Vert_{3/2}< \infty \, .\eeq We
first notice that
$$
\bkE  ( \bkE (S_{n - p_n}|{\mathcal F_{n-p_n}})\bkE (S_{2n}-S_n
|{\mathcal F_n}) | {{\mathcal F_0}}  ) = \bkE  ( \bkE (S_{n -
p_n}|{\mathcal F_{n-p_n}})\bkE (S_{2n}-S_n |{\mathcal F_{n-p_n}}) |
{{\mathcal F_0}}  ) \, ,
$$
and by stationarity and (\ref{majnormeS_n})
\begin{eqnarray*}
\Vert \bkE ( \bkE (S_{n - p_n}|{\mathcal F_{n-p_n}})\bkE (S_{2n}-S_n
|{\mathcal F_{n-p_n}}) | {{\mathcal F_0}} ) \Vert_{3/2} & \leq &
\Vert S_{n - p_n} \Vert_{3} \Vert \bkE (S_{2n}-S_n |{\mathcal
F_{n-p_n}}) \Vert_{3} \\& \leq & C \sqrt{n} \Vert \bkE
(S_{n+p_n}-S_{p_n} |{\mathcal F_{0}}) \Vert_{3}\, .
\end{eqnarray*}
Hence (\ref{resteadapSpn2}) will hold provided that \beq
\label{Condcobp3adapracinen} \sum_{n \geq 1} \frac 1n \Big \Vert
\sum_{k \geq [\sqrt n] }  \bkE (X_{k} |{\mathcal F_0})  \Big \Vert_3
< \infty\, .\eeq The fact that (\ref{Condcobp3adapracinen}) holds
under the first part of the condition (\ref{Condcobp3adap}) follows
from the following elementary lemma applied to $h(x)=\Vert \sum_{k
\geq [x] } \bkE (X_{k} |{\mathcal F_0}) \Vert_3$.
\begin{Lemma} \label{teclmahn}
Assume that $h$ is a positive function on ${\mathbb R}^+$
satisfying $ h(\sqrt{x+1}) = h(\sqrt n)$ for any $x$ in $[n-1,n[$.
Then $\sum_{n \geq 1} n^{-1} h( \sqrt n) < \infty$ if and only if
$\sum_{n \geq 1} n^{-1} h(n) <  \infty$.
\end{Lemma}
It remains to show that
\beq
\label{cond2doubleproduit} \sum_{n=1}^{\infty}\frac{1}{n^{3/2}}\Vert \bkE (
S_n \Tilde R_n  |{\mathcal F}_0) \Vert_{3/2} < \infty \, .
\eeq
Write
\begin{eqnarray*}
S_n \tilde R_n  & = & S_n \Big( \sum_{k \geq 0 } (X_{-k}-\bkE
(X_{-k} |{\mathcal F_0})) - \sum_{k \geq -n }
(X_{-k}-\bkE (X_{-k} |{\mathcal F_n})) \Big)\\
& = & S_n \Big( \bkE (S_n |{\mathcal F_n}) - S_n + \sum_{k \geq 0
} (\bkE (X_{-k} |{\mathcal F_n}) - \bkE (X_{-k} |{\mathcal F_0}))
\Big )\, .
\end{eqnarray*}
Notice first that
\begin{eqnarray*}
\Vert \bkE \big ( S_n( S_n  -\bkE (S_n|{\mathcal F_n})) |{\mathcal
F}_0 \big ) \Vert_{3/2} &  = & \Vert \bkE \big ( ( S_n -\bkE
(S_n|{\mathcal F_n}))^2 |{\mathcal F}_0 \big ) \Vert_{3/2}
 \\
&  \leq &   \Vert S_n  -\bkE (S_n|{\mathcal F_n}) \Vert_3^2  \, ,
\end{eqnarray*}
which is bounded under the second part of the condition
(\ref{Condcobp3adap}).  Now for $p_n=[\sqrt n]$, we write
$$
\sum_{k \geq 0 } (\bkE (X_{-k} |{\mathcal F_n}) - \bkE (X_{-k}
|{\mathcal F_0})) = \sum_{k \geq 0 } (\bkE (X_{-k} |{\mathcal F_n})
- \bkE (X_{-k} |{\mathcal F_{p_n}})) + \sum_{k \geq 0 } (\bkE
(X_{-k} |{\mathcal F_{p_n}}) - \bkE (X_{-k} |{\mathcal F_{0}})).
$$
Note that \begin{eqnarray*}\Big \Vert \sum_{k \geq 0 } (\bkE (X_{-k}
|{\mathcal F_{p_n}}) - \bkE (X_{-k} |{\mathcal F_{0}})) \Big \Vert_3
&  = &  \Big \Vert \sum_{k \geq 0 } (X_{-k} - \bkE (X_{-k}
|{\mathcal F_{0}})) - \sum_{k \geq 0 } (X_{-k} - (\bkE (X_{-k}
|{\mathcal
F_{p_n}}) ) \Big \Vert_3 \\
& \leq & \Big\Vert \sum_{k \geq 0 } (X_{-k} -  \bkE (X_{-k}
|{\mathcal F_{0}}))  \Big \Vert_3  + \Big \Vert  \sum_{k \geq p_n }
(X_{-k} - (\bkE (X_{-k} |{\mathcal F_{0}}) ) \Big \Vert_3 \, ,
\end{eqnarray*} which is bounded
under the second part of the condition (\ref{Condcobp3adap}). Next,
since the random variable $\sum_{k \geq 0 } (\bkE (X_{-k} |{\mathcal
F_{p_n}}) - \bkE (X_{-k} |{\mathcal F_{0}}))$ is ${\mathcal
F_{p_n}}$-measurable, we get
\begin{eqnarray*}&& \Big \Vert \bkE \Big ( S_n \sum_{k \geq 0 } (\bkE (X_{-k}
|{\mathcal F_{p_n}}) - \bkE (X_{-k} |{\mathcal F_{0}})) |{\mathcal
F_{0}} \Big ) \Big \Vert_{3/2} \\&  &  \leq \Big \Vert \bkE \Big (
S_{p_n} \sum_{k \geq 0 } (\bkE (X_{-k} |{\mathcal F_{p_n}}) - \bkE
(X_{-k}
|{\mathcal F_{0}}))  |{\mathcal F_{0}} \Big ) \Big \Vert_{3/2} \\
& & \quad \quad +  \Vert \bkE ( S_n -S_{p_n}|{\mathcal F_{p_n}})
\Vert_3 \Big \Vert \sum_{k \geq 0 } (\bkE (X_{-k}
|{\mathcal F_{p_n}}) - \bkE (X_{-k} |{\mathcal F_{0}})) \Big \Vert_3 \\
&  & \leq \Big ( \Vert  S_{p_n}\Vert_3 +  \Vert \bkE ( S_{n
-p_n}|{\mathcal F_{0}}) \Vert_3 \Big )  \Big\Vert \sum_{k \geq 0 }
(\bkE (X_{-k} |{\mathcal F_{p_n}}) - \bkE (X_{-k} |{\mathcal
F_{0}}))
\Big \Vert_3 \leq C \sqrt{p_n}  \, ,
\end{eqnarray*}
by using (\ref{Condcobp3adap}) and (\ref{majnormeS_n}). Hence, since
$p_n = [\sqrt n ]$, we get that
$$
\sum_{n=1}^{\infty}\frac{1}{n^{3/2}}  \Big \Vert \bkE \Big ( S_n
\sum_{k \geq 0 } (\bkE (X_{-k} |{\mathcal F_{p_n}}) - \bkE (X_{-k}
|{\mathcal F_{0}}))  \Big |{\mathcal F_{0}} \Big ) \Big \Vert_{3/2}
< \infty \, .
$$
It remains to show that \beq \label{proofp3nonad}
\sum_{n=1}^{\infty}\frac{1}{n^{3/2}} \Big \Vert \bkE \Big ( S_n
\sum_{k \geq 0 } (\bkE (X_{-k} |{\mathcal F_n}) - \bkE (X_{-k}
|{\mathcal F_{p_n}}))  \Big |{\mathcal F_{0}} \Big )  \Big
\Vert_{3/2} < \infty \, . \eeq  Note first that
\begin{eqnarray*}\Big \Vert \sum_{k \geq 0 } (\bkE (X_{-k} |{\mathcal
F_n}) - \bkE (X_{-k} |{\mathcal F_{p_n}})) \Big \Vert_3 &  = & \Big
\Vert \sum_{k \geq 0 } (X_{-k} -  \bkE (X_{-k} |{\mathcal F_{n}})) -
\sum_{k \geq 0 } (X_{-k} - (\bkE (X_{-k} |{\mathcal
F_{p_n}}) ) \Big \Vert_3 \\
& \leq & \Big \Vert \sum_{k \geq n } (X_{-k} -  \bkE (X_{-k}
|{\mathcal F_{0}}))  \Vert_3  + \Vert  \sum_{k \geq p_n } (X_{-k} -
(\bkE (X_{-k} |{\mathcal F_{0}}) ) \Big \Vert_3 \, .
\end{eqnarray*}
It follows that
\begin{multline*}
\Big \|\bkE \Big ( S_n \sum_{k \geq 0 } (\bkE (X_{-k} |{\mathcal
F_n}) - \bkE (X_{-k} |{\mathcal F_{p_n}}))  |{\mathcal F_{0}} \Big )
\Big \Vert_{3/2} \\\leq  C \sqrt n \Big (\Big \vert  \sum_{k \geq
p_n } (X_{-k} - (\bkE (X_{-k} |{\mathcal F_{0}}) )+  \Big \Vert_3 +
\Big \Vert  \sum_{k \geq n } (X_{-k} - (\bkE (X_{-k} |{\mathcal
F_{0}}) ) \Big \Vert_3\Big )\, .
\end{multline*}
by taking into account (\ref{majnormeS_n}). Consequently
(\ref{proofp3nonad}) will follow as soon as $$ \sum_{n \geq 1} \frac
1n \Big \Vert \sum_{k \geq [\sqrt n] } ( X_{-k} - \bkE (X_{-k}
|{\mathcal F_0}) ) \Big \Vert_3  < \infty\, ,$$ which holds under
the second part of the condition (\ref{Condcobp3adap}), by applying
Lemma \ref{teclmahn} with $h(x) = \Vert \sum_{k \geq [x] } ( X_{-k}
- \bkE (X_{-k} |{\mathcal F_0}) ) \Vert_3 $. This ends the proof of
the theorem.

\section{Appendix} \label{SectionAppendix} \setcounter{equation}{0}
\subsection{A smoothing lemma. }
\par\ssk
\begin{Lemma}\label{regul} Let $r>0$ and  $f$ be a function such
that $|f|_{\Lambda_r}< \infty$ (see Notation \ref{seminorm} for the
definition of the seminorm $|\cdot|_{\Lambda_r}$). Let $\phi_t$ be
the density of the law $N(0, t^2)$. For any real $p \geq r$ and any
positive $t$, $|f * \phi_t|_{\Lambda_p} \leq c_{r,p} t^{r-p}
|f|_{\Lambda_r}$ for some positive constant $c_{r,p}$ depending only
on $r$ and $p$. Furthermore $c_{r,r}=1$.
\end{Lemma}
\begin{Remark}
In the case where $p$ is a positive integer, the result of Lemma
\ref{regul} can be written as $\|f * \phi_t^{(p)}\|_{\infty} \leq
c_{r,p} t^{r-p} |f|_{\Lambda_r} $ .

\end{Remark}

\par\msk\rm\no{\bf Proof of Lemma \ref{regul}.}
Let $j$ be the integer such that $j<r\leq j+1$. In the case where
$p$ is a positive integer, we have
$$(f * \phi_t)^{(p)} (x) =\int \big(  f^{(j)}(u) - f^{(j)}(x) \big )
\phi_t^{(p-j)}(x-u) du  \quad \quad \text{ since $p-j\geq 1$} \, .
$$
Since $|f^{(j)}(u) - f^{(j)}(x)|\leq |x-u|^{r-j}|f|_{\Lambda_r}$, we
obtain that
$$
|(f * \phi_t)^{(p)} (x) | \leq |f|_{\Lambda_r} \int |x-u|^{r-j}
|\phi_t^{(p-j)}(x-u) |du \leq  |f|_{\Lambda_r} \int |u|^{r-j}
|\phi_t^{(p-j)}(u) |du \, .
$$
Using that $\phi_t^{(p-j)} (x) = t^{-p+j-1} \phi_1^{(p-j)} (x/t)$,
we conclude that Lemma \ref{regul} holds with the constant
$c_{r,p}=\int|z|^{r-j}\phi_1^{p-j}(z) dz $.

The case $p=r$ is straightforward. In the case where $p$ is such
that $j < r < p < j+1$, by definition
$$
|f^{(j)} * \phi_t (x) - f^{(j)} * \phi_t (y) | \leq | f
|_{\Lambda_r} |x-y|^{r-j} \, .
$$
Also, by Lemma \ref{regul} applied with $p=j+1$,
$$
|f^{(j)} * \phi_t (x) - f^{(j)} * \phi_t (y) | \leq |x-y| \|
f^{(j+1)} * \phi_t\|_{\infty}\leq
|f|_{\Lambda_r}c_{r,j+1}t^{r-j-1}|x-y| \, .
$$
Hence by interpolation,
$$
|f^{(j)} * \phi_t (x) - f^{(j)} * \phi_t (y) | \leq
|f|_{\Lambda_r}t^{r-p}c_{r,j+1}^{(p-r)/(j+1-r)}|x-y|^{p-j} \, .
$$
It remains to consider the case where $r \leq i < p \leq i+1$. By
Lemma  \ref{regul} applied successively with $p=i$ and $p=i+1$, we
obtain that
$$
|f^{(i+1)} * \phi_t (x)| \leq |f|_{\Lambda_r}c_{r,i+1}t^{r-i-1} \
\text{ and } \ |f^{(i)} * \phi_t (x)| \leq
|f|_{\Lambda_r}c_{r,i}t^{r-i} \, .
$$
Consequently
$$
|f^{(i)} * \phi_t (x) - f^{(i)} * \phi_t (y) | \leq | f
|_{\Lambda_r} t^{r-i}( 2 c_{r,i} \wedge c_{r,i+1}t^{-1} |x-y | )\, ,
$$ and by interpolation,
$$
|f^{(i)} * \phi_t (x) - f^{(i)} * \phi_t (y) | \leq | f
|_{\Lambda_r}t^{r-p}(2c_{r,i})^{1-p+i}c_{r, i+1}^{p-i}|x-y|^{p-i}\,
.
$$

\subsection{Covariance inequalities.}

In this section, we give an upper bound for the expectation of the
product of $k$ centered random variables $\Pi_{i=1}^k (X_i-{\mathbb
E}(X_i))$.
\begin{Proposition}\label{propineq}
Let $X=(X_1, \cdots, X_k)$ be a random variable with values in
${\mathbb R}^k$. Define  the number
\begin{eqnarray}\label{phi}
\phi^{(i)}& = & \phi(\sigma(X_i), X_1, \ldots, X_{i-1},X_{i+1},
\ldots, X_k)  \\ &=& \sup_{x^{} \in {\mathbb R}^{k} }\Big \|{\mathbb
E}\Big(\prod_{j=1, j\neq i}^k (\I_{X_j
> x_i} -{\mathbb P}(X_i>x_i)) | \sigma(X_i)\Big) - {\mathbb E}\Big(\prod_{j=1, j\neq i}^k (\I_{X_j
> x_i} -{\mathbb P}(X_i>x_i))  \Big ) \Big \|_{\infty} . \nonumber
\end{eqnarray}
Let $F_i$ be the distribution function of $X_i$ and $Q_i$ be the
quantile function of $|X_i|$ (see Section \ref{MC} for the
definition). Let $F_i^{-1}$ be the generalized inverse of $F_i$ and
let $D_i(u)=(F_i^{-1}(1-u)-F_i^{-1}(u))_+$. We have the inequalities
\begin{equation}\label{ineq1}
   \Big|{\mathbb E}\Big(\prod_{i=1}^k X_i-{\mathbb E}(X_i)\Big)\Big| \leq \int_0^{1}
\Big(\prod_{i=1}^k D_i(u/\phi^{(i)}) \Big) du  \,
\end{equation}
and
\begin{equation}\label{ineq2}
   \Big|{\mathbb E}\Big(\prod_{i=1}^k X_i-{\mathbb E}(X_i)\Big)\Big| \leq 2^k \int_0^{1}
\Big(\prod_{i=1}^k Q_{i}(u/\phi^{(i)}) \Big) du  \, .
\end{equation}
In addition, for any $k$-tuple $(p_1, \ldots, p_k)$ such that $1/p_1
+ \ldots + 1/p_k = 1$, we have
\begin{equation}\label{ineq2hold}
   \Big|{\mathbb E}\Big(\prod_{i=1}^k X_i-{\mathbb E}(X_i)\Big)\Big| \leq 2^k
   \prod_{i=1}^k (\phi^{(i)})^{1/p_i}\|X_i \|_{p_i}
 \, .
\end{equation}
\end{Proposition}

\noindent{\bf Proof of Proposition \ref{propineq}.} We have that
\begin{equation}\label{Hoeff}
  {\mathbb E}\Big(\prod_{i=1}^k X_i-{\mathbb E}(X_i)\Big)=
\int {\mathbb E}\Big(\prod_{i=1}^k \I_{X_i>x_i}-{\mathbb
P}(X_i>x_i)\Big) dx_1 \ldots dx_k \, .
\end{equation}
Now for all $i$,
\begin{eqnarray*}
& & {\mathbb E}\Big(\prod_{i=1}^k \I_{X_i>x_i}-{\mathbb
P}(X_i>x_i)\Big) \\
& & = {\mathbb E} \left( \I_{X_i>x_i} \Big( {\mathbb
E}\Big(\prod_{j=1, j\neq i}^k (\I_{X_j
> x_i} -{\mathbb P}(X_i>x_i)) | \sigma(X_i)\Big) - {\mathbb E}\Big(\prod_{j=1, j\neq i}^k (\I_{X_j
> x_i} -{\mathbb P}(X_i>x_i)) \Big ) \Big ) \right) \\
& & = {\mathbb E} \left( \I_{X_i \leq x_i} \Big( {\mathbb
E}\Big(\prod_{j=1, j\neq i}^k (\I_{X_j
> x_i} -{\mathbb P}(X_i>x_i)) | \sigma(X_i)\Big) - {\mathbb E}\Big(\prod_{j=1, j\neq i}^k (\I_{X_j
> x_i} -{\mathbb P}(X_i>x_i)) \Big ) \Big ) \right) .
\end{eqnarray*}
Consequently, for all $i$,
\begin{equation}\label{ineq3}
{\mathbb E}\Big(\prod_{i=1}^k \I_{X_i>x_i}-{\mathbb
P}(X_i>x_i)\Big) \leq \phi^{(i)} {\mathbb P} ( X_i \leq x_i)
\wedge {\mathbb P} ( X_i > x_i) \, .
\end{equation}
Hence, we obtain from (\ref{Hoeff}) and (\ref{ineq3}) that
\begin{eqnarray*}
 \Big|{\mathbb E}\Big(\prod_{i=1}^k X_i-{\mathbb E}(X_i)\Big)\Big| &\leq &
 \int_0^{1}
 \Big(\prod_{i=1}^k \int \I_{u/\phi^{(i)} < {\mathbb P}(X_i>x_i)}\I_{u/\phi^{(i)} \leq {\mathbb P}(X_i \leq x_i)} dx_i\Big) du \\
& \leq &  \int_0^{1}
 \Big(\prod_{i=1}^k \int \I_{F_i^{-1}(u/\phi^{(i)}) \leq x_i < F^{-1}_i(1-u/\phi^{(i)})} dx_i\Big) du ,
\end{eqnarray*}
and (\ref{ineq1}) follows. Now (\ref{ineq2}) comes from
(\ref{ineq1}) and the fact that $D_i(u) \leq 2 Q_{i} (u)$ (see Lemma
6.1 in Dedecker and Rio (2006)). $\square$

\begin{Definition} For a quantile function $Q$ in ${\mathbb
L}_1([0,1] , \lambda )$, let ${\mathcal F}(Q, P_X)$ be the set of
functions $f$ which are nondecreasing on some open interval of
${\mathbb R}$ and null elsewhere and such that $Q_{|f(X)|} \leq
Q$.\\ Let ${\mathcal C}(Q, P_X)$ denote the set of convex
combinations $\sum_{i=1}^{\infty} \lambda_i f_i$ of functions $f_i$
in ${\mathcal F}(Q, P_X)$ where $\sum_{i=1}^{\infty} |\lambda_i|
\leq 1$ (note that the series $\sum_{i=1}^{\infty} \lambda_i f_i(X)$
converges almost surely and in ${\mathbb L}_1(P_X)$).
\end{Definition}

\begin{Corollary}\label{propineq2}
Let $X=(X_1, \cdots, X_k)$ be a random variable with values in
${\mathbb R}^k$ and let the $\phi^{(i)}$'s be defined by
(\ref{phi}). Let $(f_i)_{1\leq i \leq k}$ be $k$ functions from
${\mathbb R}$ to ${\mathbb R}$, such that $f_i \in {\mathcal C}(Q_i,
P_{X_i})$.  We have the inequality
$$
   \Big|{\mathbb E}\Big(\prod_{i=1}^k f_i(X_i)-{\mathbb E}(f_i(X_i))\Big)\Big|
\leq 2^{2k-1} \int_0^{1} \prod_{i=1}^k Q_{i}\Big
(\frac{u}{\phi^{(i)}} \Big ) du \, .
$$
\end{Corollary}

\noindent{\bf Proof of Corollary \ref{propineq2}.} Write for all $1
\leq i \leq k$, $f_i = \sum_{j = 1}^{\infty} \lambda_{j, i} f_{j
,i}$ where $\sum_{j=1}^{\infty} |\lambda_{j,i}| \leq 1$ and $f_{j,i}
\in {\mathcal F}(Q_i, P_{X_i})$. Clearly
\begin{eqnarray}\label{clear}
\Big|{\mathbb E}\Big(\prod_{i=1}^k f_i(X_i)-{\mathbb
E}(f_i(X_i))\Big)\Big| & \leq & \sum_{j_1=1}^{\infty}\cdots
\sum_{j_k=1}^\infty \Big(\prod_{i=1}^k|\lambda_{j_i,
i}|\Big)\Big|{\mathbb E}\Big(\prod_{i=1}^k f_{j_i ,
i}(X_i)-{\mathbb E}(f_{j_i , i}(X_i))\Big)\Big| \nonumber \\
& \leq & \sup_{j_1 \geq 1, \ldots, j_k \geq 1}\Big|{\mathbb
E}\Big(\prod_{i=1}^k f_{j_i , i}(X_i)-{\mathbb E}(f_{j_i ,
i}(X_i))\Big)\Big|\,  .
\end{eqnarray}
Since each $f_{j_i , i }$ is nondecreasing on some interval,
$$
\phi(\sigma(f_{j_i , i}(X_i)), f_{j_1 , 1}(X_1), \ldots, f_{j_{i-1}
, i-1}(X_{i-1}),f_{j_{i+1} , i+1}(X_{i+1}), \ldots, f_{j_{k} ,
k}(X_{k})) \leq 2^{k-1} \phi^{(i)}\, .
$$ Then applying (\ref{ineq2}) on the right hand side of (\ref{clear}), we derive that
$$
   \Big|{\mathbb E}\Big(\prod_{i=1}^k f_i(X_i)-{\mathbb E}(f_i(X_i))\Big)\Big|
\leq 2^{k} \int_0^{1} \prod_{i=1}^k Q_{i}\Big
(\frac{u}{2^{k-1}\phi^{(i)}} \Big ) du \, ,
$$ and the result follows by a change-of-variables. $\square$

\medskip

Recall that  for any $p \geq 1$, the class ${\mathcal C} (p, M, P_X
)$ has been introduced in  the definition \ref{defclosedenv}.

\begin{Corollary}\label{propineq3} Let $X=(X_1, \cdots, X_k)$ be a random variable with values in
${\mathbb R}^k$ and let the $\phi^{(i)}$'s be defined by
(\ref{phi}). Let a $k$-tuple $(p_1, \ldots, p_k)$ such that $1/p_1 +
\ldots + 1/p_k = 1$ and let $(f_i)_{1\leq i \leq k}$ be $k$
functions from ${\mathbb R}$ to ${\mathbb R}$, such that $f_i \in
{\mathcal C}(p_i, M_i, P_{X_i})$.  We have the inequality
$$
   \Big|{\mathbb E}\Big(\prod_{i=1}^k f_i(X_i)-{\mathbb E}(f_i(X_i))\Big)\Big|
\leq 2^{2k-1}  \prod_{i=1}^k (\phi^{(i)})^{1/p_i} M_i^{1/p_i}\, .
$$
\end{Corollary}


\end{document}